\documentclass[a4paper,11pt,reqno]{amsart}
\usepackage{amsmath}
\usepackage{amssymb}
\usepackage{amsfonts}
\usepackage{graphicx}
\usepackage{mathtools}
\usepackage[colorlinks]{hyperref}

\renewcommand\eqref[1]{(\ref{#1})} 
\graphicspath{ {images/} }
\usepackage{thmtools}

\newtheorem{theorema}{Theorem}[section]

\title[Sharp Hardy Inequality and Sobolev--Lorentz Embeddings] {Sharp Quantitative Forms of the Hardy Inequality on Cartan--Hadamard Manifolds via Sobolev--Lorentz Embeddings}

\author[A. Banerjee]{Avas Banerjee}
	\address[A. Banerjee]{Theoretical Statistics and Mathematics Unit, Indian Statistical Institute, Delhi Center, S.J. Sansanwal Marg, New Delhi, Delhi 110016, India
	}
	\email{avas24r@isid.ac.in}
    
\author[D. Ganguly]{Debdip Ganguly}
\address[D. Ganguly]{Theoretical Statistics and Mathematics Unit, Indian Statistical Institute, Delhi Center, S.J. Sansanwal Marg, New Delhi, Delhi 110016, India
    }
    \email{debdip@isid.ac.in}
    
\author[P. Roychowdhury]{Prasun Roychowdhury}
    \address[P. Roychowdhury]{Department of Mathematics, Indian Institute of Technology Hyderabad, Kandi, Sangareddy, Telangana, 502285, India.
    } 
	\email{prasunrc@math.iith.ac.in}


\date{\today}

\theoremstyle{plain}
\newtheorem{theorem}{Theorem}[section]
\newtheorem{proposition}{Proposition}[section]
\newtheorem{lemma}{Lemma}[section]

\newtheorem{remark}{Remark}[section]
\newtheorem{definition}{Definition}[section]

\numberwithin{equation}{section} \allowdisplaybreaks

\usepackage[text={6.8in,10in},centering]{geometry}
\newcommand{\hn}{\mathbb{H}^N}
\newcommand{\dv}{\: {\rm d}v_{\hn}}

\newcommand{\dvg}{\: {\rm d}v_{g}}
\newcommand{\gradg}{\nabla_g}

\newcommand{\lapg}{\Delta_{g}}

\newcommand{\m}{\mathbb{M}^N}
\newcommand{\dx}{\:{\rm d}x}
\newcommand{\dy}{\:{\rm d}y}
\newcommand{\dt}{\:{\rm d}t}
\newcommand{\ds}{\:{\rm d}s}
\newcommand{\dr}{\:{\rm d}r}
\newcommand{\dm}{\:{\rm d}m}
\newcommand{\dsn}{\:{\rm d}\theta}
\newcommand{\sn}{\mathbb{S}^{N-1}}
\newcommand{\rn}{\mathbb{R}^N}
\newcommand{\uv}{\omega_N}
\newcommand{\drho}{\:{\rm d}\varrho}
\newcommand{\si}{\mathcal{I}}
\newcommand{\st}{\mathcal{T}}
\newcommand{\scal}{\mathcal{S}}
\newcommand{\eu}{\mathcal{E}}
\newcommand{\sd}{\mathcal{D}}

\begin{document}
\begin{abstract}
In this article, we investigate the quantitative form of the classical Hardy inequality. In our first result, we prove the following quantitative bound under the assumption that the $\m$ is a Riemannian model satisfying the \emph{Pólya--Szegő inequality}: Specifically, we prove that
\[
\|\nabla_g u\|^2_{L^{2}(\m)} -  \frac{(N-2)^2}{4}\left\|\frac{u}{r(x)}\right\|^2_{L^2(\m)} \geq C \; [\mbox{dist}(u, \mathcal{Z})]^{\frac{4N}{N-2}}\left\|\frac{u}{r(x)}\right\|^2_{L^2(\m)},
\] 
for every real-valued weakly differentiable function $u$ on $\m$ such that $|\nabla_g u| \in L^2(\m)$ and $u$
decays to zero at infinity. Here $r(x) = d_g(x, x_0)$ denotes the geodesic distance from a fixed pole $x_0,$  the set $\mathcal{Z}$ represents the family of ``virtual" extremals, and the distance is understood in an appropriate generalized Lorentz-type space. Our approach is built on the symmetrization technique on manifolds, combined with a novel Jacobian-type transformation that provides a precise way for comparing volume growth, level sets, and gradient terms across the two geometries of Euclidean and manifold settings, respectively. When coupled with symmetrization, this framework yields sharp control over the relevant functionals and reveals how the underlying curvature influences extremal behaviour. Our result generalizes the seminal result of Cianchi-Ferone [Ann. Inst. H. Poincar\'e C Anal. Non Lin\'eaire 25 (2008)] to the curved spaces. Moreover, building upon this transformation, we introduce a Lorentz-type framework naturally adapted to the volume growth of Cartan--Hadamard model manifolds and establish an \emph{optimal} Sobolev--Lorentz embedding in this setting, including the identification of the corresponding sharp constant and its non-attainment. Finally, we establish a quantitative correspondence between the Hardy \emph{deficit} on the manifold and an appropriate weighted Hardy deficit in Euclidean space, showing that each controls the other. In particular, this provides a quantitative link between the Hardy inequality on the curved manifold and its Euclidean counterpart.
\end{abstract}
 
 \maketitle

\section{Introduction}
In 1920, in his study of Fourier series, \textsc{G.~H.~Hardy} introduced the following integral inequality, and a complete proof appeared five years later. Let $f$ be a nonnegative measurable function on $\mathbb{R}_+$, and define 
\[
    F(x) = \int_0^x f(t)\,\mathrm{d}t.
\]
Then, the sharp Hardy inequality reads 
\begin{equation}\label{Hardy-int}
    \int_0^\infty \left( \frac{F(x)}{x} \right)^2 \,\mathrm{d}x 
    \;\le\; 4 \int_0^\infty |f(x)|^2 \,\mathrm{d}x,
\end{equation}
whenever the right-hand side is finite; see \cite[Theorem~327]{HPL} for background and further discussion.

A simple computation shows that \eqref{Hardy-int} is equivalent to the one-dimensional differential form of Hardy's inequality,
\begin{equation}\label{Hardy-diff}
    \int_0^\infty \left( \frac{f(x)}{x} \right)^2 \mathrm{d}x
    \;\le\; 4 \int_0^\infty |f'(x)|^2 \mathrm{d}x
\end{equation}
which holds for sufficiently regular functions vanishing at infinity. 

Hardy-type inequalities such as \eqref{Hardy-int}--\eqref{Hardy-diff} play a central role in analysis and partial differential equations, serving as prototypes for a wide class of functional inequalities. Over the past century, many refinements, extensions, and geometric generalizations of Hardy's original inequality have been developed. In particular, the most natural higher-dimensional analogue asserts that if $N \geq 3,$  then 
\begin{equation}\label{intro-hardy}
\frac{(N-2)^2}{4}\int_{\mathbb{R}^N} \dfrac{|u|^2}{|x|^2} \; {\rm d}x \; \leq \; \int_{\mathbb{R}^N} |\nabla u|^2 \; {\rm d}x,
\end{equation}
for every real-valued weakly differentiable function $u \in D^{1, 2}(\mathbb{R}^N),$ where 
\[
D^{1, 2}(\mathbb{R}^N) : = \{ u : \mathbb{R}^N \rightarrow \mathbb{R} : |\nabla u| \in L^2(\mathbb{R}^N), 
\mbox{and} \ \mathcal{L}^N ( \{ u \geq t \} ) < \infty, \ \mbox{for every} \ t >0 \}.
\]
Here $\mathcal{L}^N(A)$ denotes the Lebesgue measure of a set $A \subset \mathbb{R}^N.$ Equivalently, $D^{1, 2}(\mathbb{R}^N)$ is the completion of $C_c^{\infty}(\mathbb{R}^N)$ with respect to norm 
\[
\| u\| = \left( \int_{\mathbb{R}^N} |\nabla u|^2 \; {\rm d}x \right)^{\frac{1}{2}}.
\]
Moreover, define the Rayleigh quotient
\[
R[u]=\frac{\displaystyle\int_{\mathbb{R}^N}|\nabla u|^2\,{\rm d}x}
{\displaystyle\int_{\mathbb{R}^N}\frac{|u|^2}{|x|^2}\,{\rm d}x},
\qquad u\in D^{1,2}(\mathbb{R}^N)\setminus\{0\}.
\]
The optimal Hardy constant is the infimum, $C_{\rm H}=\inf_{u\in D^{1,2}(\mathbb{R}^N)\setminus\{0\}} R[u],
$ and it is well known that  $C_{\rm H} = \frac{(N-2)^2}{4}.$  To this end, the formal Euler--Lagrange equation for minimizers of \(R[u]\) is
\[
-\Delta u = \lambda\frac{u}{|x|^2}
\]
for some Lagrange multiplier \(\lambda\). The radial fundamental solution of the homogeneous ODE associated with the above is proportional to
\[
u(x)\sim |x|^{-\frac{N-2}{2}}.
\]
Consequently, a function that would make the equality in \eqref{intro-hardy} must behave like \(|x|^{-(N-2)/2}\) near the origin. But \(|x|^{-(N-2)/2}\notin L^{2^*}_{\rm loc}\) and, more relevantly, does not belong to \(D^{1,2}(\mathbb{R}^N)\) (its gradient is not square integrable at \(0\) or at infinity). Hence, the infimum is not achieved. Because equality is not achieved, one can often strengthen \eqref{intro-hardy} by adding a positive remainder term on the left. Rewrite \eqref{intro-hardy} as a nonnegativity statement for the Schr\"odinger operator with inverse-square potential:
\[
\int_{\mathbb{R}^N}|\nabla u|^2\,{\rm d}x - c\int_{\mathbb{R}^N}\frac{|u|^2}{|x|^2}\,{\rm d}x \ge 0
\quad\text{for all }u\in C_c^\infty(\mathbb{R}^N)
\]
if and only if \(c\le C_{\rm H}=(N-2)^2/4\). Thus \(c=(N-2)^2/4\) is the \emph{critical} coupling constant of the potential \(-c|x|^{-2}\). For \(c>C_{\rm H}\) the quadratic form becomes indefinite and the operator \( -\Delta - c|x|^{-2}\) admits negative spectrum (or is unbounded below). For \(c=C_{\rm H}\), the operator is at the threshold of nonnegativity; this criticality is reflected in the non-attainment of the infimum and in the existence of singular null-sequences. From the above discussion, it is straightforward to verify that the functions
\[
\mathcal{V}_a(x) = a\, |x|^{\frac{2-N}{2}}, \qquad x \in \mathbb{R}^N \setminus \{0\},
\]
with \( a \in \mathbb{R} \setminus \{0\} \) and \( N \ge 3 \), satisfy the Euler--Lagrange equation
\[
-\Delta u = \frac{(N-2)^2}{4}\,\frac{u}{|x|^2}, \qquad x \in \mathbb{R}^N \setminus \{0\}.
\]
Hence, the family \( \{\mathcal{V}_a\}_{a \ne 0} \) represents the natural candidates for extremals of~\eqref{intro-hardy}. However, since these functions do not belong to \( D^{1,2}(\mathbb{R}^N) \), the inequality is never attained within this space. Following the terminology introduced in the literature, we refer to such functions as \emph{virtual extremals}. 
\medskip

{\bf Question.} A natural question that arises is whether these virtual extremals belong to some natural function space.

To address the above question, we first recall a classical result due to Peetre, which provides an improvement over the Sobolev inequality: 

\begin{align*}
  D^{1,2}(\mathbb{R}^N)
  \hookrightarrow 
  L^{2^\star, 2}(\mathbb{R}^N)
  \hookrightarrow 
  L^{2^\star}(\mathbb{R}^N),
\end{align*}
where the first embedding is known as \emph{Peetre’s embedding}, which refines the classical Sobolev inequality and  $ L^{2^\star, 2}(\mathbb{R}^N)$ defines the Lorentz space.  In the seminal paper, Cianchi and Ferone \cite{cf-aihp} reformulated the Hardy inequality in terms of the Lorentz norm. In particular, via the symmetrization argument, inequality \eqref{intro-hardy} is equivalent to the Lorentz-Norm inequality 
$$
\alpha_N^{1/N} \dfrac{(N-2)}{2} \| u \|_{L^{2^\star, 2}(\mathbb{R}^N)} \leq \|\nabla u\|_{L^{2}(\mathbb{R}^N)},
$$
where $\alpha_N$ is the Euclidean measure of the unit ball. In view of the above inequality, one may naturally ask the following stability question: the norm $\|\cdot\|_{L^{2^\star, 2}(\mathbb{R}^N)}$ can be regarded as the natural choice to measure the distance from the extremizer in the following sense:
$$
\|\nabla u\|_{L^{2}(\mathbb{R}^N)} \;- \;\alpha_N^{1/N} \dfrac{(N-2)}{2} \| u \|_{L^{2^\star, 2}(\mathbb{R}^N)} \geq 
C\;d(u), \quad \forall u \in D^{1, 2}(\mathbb{R}^N),
$$
where $C:= C(N)>0$ is a constant, and $d(u) := \inf_{a \in \mathbb{R}} \dfrac{\| u - \mathcal{V}_a\|_{L^{2^\star, 2}(\mathbb{R}^N)}}{\| u\|_{L^{2^\star, 2}(\mathbb{R}^N)}}.$

Unfortunately, this is not possible, since the functions $\mathcal{V}_a$, whose gradients do not belong to $L^2(\mathbb{R}^N)$, are neither contained in the Lorentz space $L^{2^\ast, 2}(\mathbb{R}^N)$ nor even in the larger Lebesgue space $L^{2^\ast}(\mathbb{R}^N)$ appearing in the classical Sobolev embedding. Indeed, $\mathcal{V}_a$ fails to satisfy the necessary integrability condition at infinity required for these spaces. This lack of integrability reflects the borderline behaviour of the extremal profile for the Hardy inequality when expressed in terms of weak-type spaces. 

In fact, the smallest rearrangement-invariant space containing $\mathcal{V}_a$ is the Marcinkiewicz space $L^{2^\ast, \infty}(\mathbb{R}^N)$, also known as the weak-$L^{2^\ast}$ space. This observation serves as an appropriate functional framework for studying refined inequalities and stability questions. In particular, we state a quantitative result due to Cianchi and Ferone in the case $p=2$. To this end, let  $D(u) := \inf_{a \in \mathbb{R}} \dfrac{\| u - \mathcal{V}_a\|_{L^{2^\star, \infty}(\mathbb{R}^N)}}{\| u\|_{L^{2^\star, 2}(\mathbb{R}^N)}}.$

\begin{theorema}[Cianchi and Ferone \cite{cf-aihp}]
Let $N \geq 3$. Then there exists a constant $C = C(N) > 0$ such that 
\begin{equation*}
\|\nabla u\|^2_{L^{2}(\mathbb{R}^N)} -  C_H\left\|\frac{u}{|x|}\right\|^2_{L^2(\mathbb{R}^N)} \geq C \; [D(u)]^{\frac{4N}{N-2}}\left\|\frac{u}{|x|}\right\|^2_{L^2(\mathbb{R}^N)},
\end{equation*}
for every real-valued weakly differentiable function $u$ on $\mathbb{R}^N$ decaying to zero at infinity and such that $|\nabla u| \in L^2(\mathbb{R}^N).$  
\end{theorema}
The proof rests on a classical symmetrization argument combined with a careful analysis of level sets and a quantitative comparison with the extremal profiles.

\medskip 

\subsection{Question} It is natural to ask whether these quantitative stability results admit an extension to Cartan--Hadamard manifolds. 
On a simply connected complete Riemannian manifold with nonpositive sectional curvature, the lack of translation invariance and the presence of curvature-dependent volume growth require a refined analysis: one must replace Euclidean symmetrization by suitable comparison arguments (or by symmetrization on model spaces), control the Jacobian factors arising in polar coordinates, and exploit geometric inequalities such as the isoperimetric-type estimates adapted to the curvature bounds. Consequently, both the form of the extremal profiles and the stability rates may depend explicitly on curvature bounds or on the manifold's asymptotic geometry.

Let $(\m, g)$ be a complete, simply connected Riemannian manifold of dimension $N\ge3$ with nonpositive sectional curvature. 
Assume further that $\m$ satisfies some volume growth condition. Then for $N \geq 3$ there exist constants $C= C(N)$ and an appropriate manifold-dependent distance $d_{\m}(u)$ to the family of extremals such that every real-valued, weakly differentiable function $u$ on $\m$ with $|\gradg u|\in L^2(\m)$ and vanishing at infinity satisfies a quantitative Hardy inequality
\[
\|\nabla_g u\|^2_{L^{2}(\m)} -  C_H\left\|\frac{u}{r(x)}\right\|^2_{L^2(\m)} \geq C \; [d_{\m}(u)]^{\frac{4N}{N-2}}\left\|\frac{u}{r(x)}\right\|^2_{L^2(\m)},
\]
where $r(x)=\mathrm{dist}_g(x,x_0)$ for a fixed reference point $x_0\in \m$. In recent years, substantial progress on the stability of related inequalities, such as the Caffarelli–Kohn–Nirenberg and Hardy–Littlewood–Sobolev inequalities, has been achieved by several authors; we only mention \cite{cfll, lu-2, CLT} without any claim of completeness and refer to the references therein.

\subsection{Riemannian Model Manifolds} In the present article, we develop a unified framework for the improvement and stability of the Hardy inequalities on manifolds which satisfy some symmetry assumption so that Euclidean symmetrization techniques transfer to the Riemannian setting most cleanly and this naturally leads to the fact that the ambient manifold enjoys rotational symmetry, so we restrict attention to model (warped-product) manifolds,  whose metric can be written in the form
\begin{equation*}
g = \dr \otimes \dr + \psi(r)^2 g_{\sn},
\end{equation*}
where $r(x)=\text{dist}(x,x_0)$ be the geodesic distance between $x\in \m$ and $x_0\in \m$, a fixed point 
(referred to as the \emph{pole} of the manifold). The vector $\dr$ 
represents the radial direction, and $g_{\sn}$ is the standard metric 
on the unit sphere $\sn$. The function 
$\psi : [0,\infty) \to [0,\infty)$ is a smooth function that fully 
determines the geometry of $\m$. Moreover, the Cartan--Hadamard assumption is equivalent to requiring that the warping function 
$\psi$ be convex. A prototypical example is given by the choice $\psi(r) = \sinh r$, 
which yields the classical model of the hyperbolic space $\mathbb{H}^N$. Another important question concerns the validity of the Pólya--Szegő inequality on model manifolds. Muratori and Volzone~\cite[Theorem~3.9]{mv} provided a class of model manifolds for which
Pólya--Szegő inequality is violated. In the same work~\cite{mv}, the authors identified 
sufficient geometric conditions ensuring the validity of the Pólya--Szegő inequality. In particular, they 
proved that for any model manifold, the P\'olya--Szegő inequality holds whenever $\m$ satisfies a \emph{centered isoperimetric inequality} (see section~\ref{sec-symm} for further details). 
\textit{Throughout the paper, our manifold $\m$ is assumed to satisfy the P\'olya--Szegő inequality.}  For a more detailed discussion, including precise definitions, geometric background, 
and the technical framework used in this work, we refer the reader to 
Section~\ref{prem}. There, we provide a systematic presentation of model manifolds, and
their isoperimetric properties that play a central role 
in our analysis. 
\subsection{Embeddings and Function spaces} Apart from the geometric framework, it is essential to develop several functional-analytic 
tools that allow us to define and work with \emph{Lorentz spaces} on such manifolds. 
In fact, the natural candidates for extremals of the Hardy inequality on these manifolds 
typically fail to belong to the standard Lorentz spaces associated with the Riemannian measure. 
This failure is due to the specific decay behavior of these extremal functions, 
which is often dictated by the volume growth of the warping function of the \emph{ model manifold}. 
To address this issue, we introduce a new \emph{Lorentz-type quasi-norm} on $\m$, 
tailored to capture the precise integrability and non-linear scaling properties of these extremals. 
These Lorentz-type spaces retain many of the classical embedding and interpolation properties, 
which allow us to extend sharp functional inequalities, study extremals, and analyze 
quantitative stability in a setting where the standard Lorentz spaces are inadequate. We 
refer to Section~\ref{sec-symm} for more details. This section also highlights how the centered isoperimetric inequality fits naturally within the broader framework of symmetrization and rearrangement methods on curved spaces.

\subsection{Highlights of our main results} For the reader's convenience, we present here the main results of the paper in a slightly formal way. These statements are not fully rigorous at this stage, as their precise 
formulation requires technical concepts introduced in Section~\ref{prem}, such as the geometric notion 
 and the formulation of function spaces in which quantitative stability type results for the Hardy inequality are intended to hold. 
For the complete and rigorous statements, we refer the reader to Section~\ref{sec-3} and Section~\ref{sec-SL}. We begin with our first stability results. 

\begin{theorem}\label{intro-main-th-hardy-cf-stab}
Consider an $N$-dimensional Cartan--Hadamard model manifold $(\m,g)$, 
with $N \ge 3$. Suppose $u\in D^{1,2}(\m)$ and denote 
     \begin{equation}\label{hardy stability distance}
        d_{\m}(u)=\inf_{a\in \mathbb{R}}\frac{\|u-\mathcal{U}_a\|_{\tilde{L}^{2^\star,\infty}(\m)}}{\|u\|_{D^{1,2}(\m)}}.
    \end{equation}
    Then there exists a constant $C=C(N)>0$, such that there holds
    \begin{equation}\label{hyperbolic hardy stability}
        \left(\frac{N-2}{2}\right)^2\int_{\m}\frac{u^2}{r^2}\dvg\left(1+C\, (d_{\m}(u))^{\frac{4N}{N-2}}\right)\leq\int_{\m}|\gradg u|^2\dvg,
    \end{equation}
    where $\tilde{L}^{2^\star,\infty}(\m)$ denotes the Lorentz-type quasi-norm (defined in \eqref{quasinorm}) and $\mathcal{U}_a$ denotes the ``virtual" extremals for the Hardy inequality in $\m$, which are explicitly given by \eqref{extremizer}
\end{theorem}

\begin{remark}
{\rm
Theorem~\ref{intro-main-th-hardy-cf-stab} provides a quantitative version of the Hardy inequality on Cartan--Hadamard model manifolds. 
The estimate~\eqref{hyperbolic hardy stability} demonstrates that the deficit in the classical Hardy inequality not only remains nonnegative but also quantitatively controls the distance 
\[
d_{\m}(u) = \inf_{a\in \mathbb{R}}\frac{\|u-\mathcal{U}_a\|_{\tilde{L}^{2^\star,\infty}(\m)}}{\|u\|_{D^{1,2}(\m)}}
\]
of a function \(u\) from the manifold’s family of “virtual” extremals \(\mathcal{U}_a\) measured in the Lorentz-type quasi-norm \(\tilde{L}^{2^\star,\infty}(\m)\). Consequently, functions that nearly saturate the Hardy inequality (formally) must lie close to these extremal profiles. 
}
\end{remark}

The analysis of \emph{optimal Sobolev--Lorentz embeddings} on model manifolds provides a refined understanding of how curvature and geometry influence the classical Sobolev framework. In contrast to the standard embedding 
\[
D^{1,2}(\m) \hookrightarrow L^{2^\ast}(\m),
\]
its Lorentz-space refinement 
\[
D^{1,2}(\m) \hookrightarrow L^{2^\ast,2}(\m) \hookrightarrow L^{2^\ast}(\m), 
\]
captures borderline integrability phenomena and yields a sharper function space. On a Cartan--Hadamard model manifold \((\m,g)\), the embedding behaviour is deeply influenced by the geometric structure in the radial volume density function \(\psi(r)\), which determines whether the embedding is continuous, compact, or optimal.  We prove an \emph{optimal} Sobolev-Lorentz embedding on the Cartan-Hadamard models. In particular, we prove

\begin{theorem}
Let $(\m,g)$ be an $N$-dimensional Cartan--Hadamard model manifold with $N \geq 3$. Then the Sobolev space $D^{1,2}(\m)$ admits the optimal Lorentz refinement
\[
D^{1,2}(\m) \hookrightarrow L^{2^\star,2}(\m) \hookrightarrow L^{2^\star}(\m),
\]
where $L^{2^\star,2}(\m)$ is the smallest rearrangement-invariant space into which $D^{1,2}(\m)$ embeds. Moreover, this embedding is continuous and satisfies the sharp inequality
\begin{align}\label{avas-best}
    \|u\|_{L^{2^\star,2}(\m)} 
    \leq S_{N,2^\star,2} \, \|u\|_{D^{1,2}(\m)},
\end{align}
where the constant $S_{N,2^\star,2}$ is given explicitly in~\eqref{alv-best}. In addition, this optimal constant is never attained by any function in $D^{1,2}(\m)$.
\end{theorem}

\medskip 

\begin{remark}
{\rm
To the best of our knowledge, the above result has been established only in the hyperbolic setting \cite{VHN}, where Nguyen \cite[Theorem~1.2]{VHN} obtained it via a different method (see also the Lorentz–Sobolev inequalities for higher-order derivatives in the hyperbolic space \cite{vhn-lorentz-poin}). However, the discussion in \cite{VHN} did not address the optimality of the embedding nor the question of attainment of the best constant. In contrast, our result not only extends this framework to general Cartan--Hadamard model manifolds but also identifies the precise Lorentz space providing the smallest rearrangement-invariant target and explicitly determines the sharp constant of the embedding.
}
\end{remark}
\medskip

Next, we establish a precise link between the geometric and Euclidean settings by showing that the Hardy inequality on a Cartan–Hadamard model manifold can be both obtained from and strengthened via a suitably weighted Euclidean Hardy inequality. Specifically, we prove that the Hardy \emph{deficit} on the manifold is quantitatively controlled from below by a weighted Hardy deficit in Euclidean space, computed for a function naturally associated with the manifold through the radial Jacobian transformation \eqref{jacobian}. As a first step in this analysis, we establish a weighted Euclidean Hardy inequality, which will serve as the fundamental tool for deriving the corresponding classical Hardy inequalities on the manifold. In other words, the Hardy deficit in the Cartan-Hadamard manifold is bounded below by a Hardy-type deficit for a certain symmetrized function in the Euclidean space. 

\begin{theorem}
Let $(\m,g)$ be a Cartan--Hadamard model manifold of dimension $N \ge 3$. 
Then, for every function $u \in D^{1,2}(\m)$, there exists a radial, decreasing 
function $F \in D^{1,2}_{\mathrm{rad}}(\mathbb{R}^N \setminus \{0\})$ such that
\begin{multline}\label{intro-manifold-eu-eq}
\int_{\m} |\nabla_g u|^2 \, \dvg
- \left(\frac{N-2}{2}\right)^2 \int_{\m} \frac{u^2}{r^2} \, \dvg \\
\ge 
\int_{\mathbb{R}^N} |\nabla F(x)|^2 \, (w(x))^{N-1} \dx
- \left(\frac{N-2}{2}\right)^2 \int_{\mathbb{R}^N} \frac{|F(x)|^2}{|x|^2} \, (w(x))^{N-1} \dx,
\end{multline}
where 
\[
w(x) = \frac{\psi(|x|)}{|x|} \quad \text{for all } x \in \mathbb{R}^N \setminus \{0\}.
\]
\end{theorem}
\begin{remark}
{\rm 
Thanks to Lemma~\ref{aux-fnl}, the function $w(x)$ appearing in Theorem~\ref{manifold-eu} 
is monotonically increasing. Consequently, by applying Proposition~\ref{wg-ecld-hyardy}, we deduce that the right-hand side of inequality~\eqref{intro-manifold-eu-eq} is nonnegative. 
This observation establishes a direct connection between the classical Hardy inequality 
on the Cartan--Hadamard model manifold and the corresponding weighted Euclidean 
Hardy inequality: in other words, the validity of the Hardy inequality~\eqref{hardy} 
on the manifold can be inferred via the weighted Euclidean framework. 

}
\end{remark}


\begin{theorem}
Let $(\m,g)$ be a Cartan--Hadamard model manifold of dimension $N \ge 3$. 
Then, for every function $u \in D^{1,2}(\mathbb{R}^N \setminus \{0\})$, there exists a 
radial, decreasing function $F$ defined on $\m \setminus \{x_0\}$, vanishing 
at infinity, such that
\begin{multline}\label{eu}
\int_{\mathbb{R}^N} |\nabla u|^2 \dx
- \left(\frac{N-2}{2}\right)^2 \int_{\mathbb{R}^N} \frac{|u|^2}{|x|^2} \dx \\[2mm]
\ge 
\int_{\m} |\nabla_g F|^2 \, \big(\Psi(r)\big)^{N-1} \, \dvg
- \left(\frac{N-2}{2}\right)^2 \int_{\m} \frac{|F|^2}{r^2} \, \big(\Psi(r)\big)^{N-1} \, \dvg,
\end{multline}
where
\begin{align*}
   \Psi(r):=\frac{r}{\psi(r)} \qquad \text{for all } r:= \mathrm{dist}_g(x,x_0)>0.
\end{align*}
\end{theorem}

\begin{remark}
{\rm
We generalize the classical Euclidean Hardy inequality to a broad class of Riemannian model manifolds and establish precise weighted inequalities that reflect the underlying curvature.
}
\end{remark}

\subsection{Strategy of proof and outline of the paper} We now provide a brief explanation of the central idea on which this work is based. Our approach is built on the symmetrization technique on the manifolds, combined with a new delicate transformation that allows us to relate geometric quantities on the manifold to their Euclidean counterparts in a constructive and analytic manner.
\begin{itemize}

\item 
For a given radius \( \tilde{r} > 0 \) in the manifold \( \m \), let \( V(\tilde{r}) \) denote the corresponding volume function as defined in~\eqref{vol}. 
We then consider the radius \( \tilde{\varrho} > 0 \) in the Euclidean space \( \mathbb{R}^N \) for which the Euclidean volume matches \( V(\tilde{r}) \). 
That is, \( \tilde{\varrho} \) is chosen so that the Euclidean ball \( B_{\tilde{\varrho}} \subset \mathbb{R}^N \) satisfies
\begin{align}\label{transformation}
  |B_{\tilde{\varrho}}|_{\mathbb{R}^N} \;=\; V(\tilde{r}).  
\end{align}
This gives a one-to-one correspondence between geodesic radii in \( \m \) and Euclidean radii determined purely by the geometry of \( \m \).

\item 
The above volume–matching procedure naturally induces a correspondence between geodesic balls in the manifold \( \m \) and Euclidean balls.  
In particular, it yields a monotone mapping \( \tilde{r} \mapsto \tilde{\varrho} \), (under certain curvature assumptions), which can be inverted to express the manifold radius in terms of its Euclidean analogue.  
This inverse relation will be central in formulating our transformation.

\item 
Under suitable volume growth assumptions on \( \m \), one can compute the Jacobian of the inverse change of variables explicitly.  
Incorporating this Jacobian allows us to define a new \emph{transformation} that correctly captures how volumes change due to the geometry of \( \m \). This transformation is delicate and remarkably effective for analytic purposes.

\item 
The resulting transformation provides a powerful bridge between the curved geometry of \( \m \) and the classical Euclidean setting.  
In particular, it enables us to translate functional inequalities on \( \m \) into their Euclidean counterparts in an exact and quantitative manner.  
Consequently, sharp Euclidean tools—such as symmetrization, comparison principles, and optimal inequalities—can be employed to study problems arising on \( \m \), with the curvature of the manifold entering solely through the Jacobian of the transformation.

\end{itemize}
	
\medskip
	
 The structure of the paper is as follows.
	
	\begin{itemize}
		
		\item[Section 1:] The introduction section contains a brief background on the Hardy inequality and highlights the main results on the quantitative stability in this geometric setting, extending the ideas of Cianchi–Ferone \cite{cf-aihp} to non-Euclidean spaces and offering new insights into the role of geometry in functional inequalities.

		\item[Section 2:] Contains the geometric preliminaries, geometric assumptions on the manifolds, and recalls some of the well-known facts concerning symmetrization technique on manifolds, e.g., Hardy-Littlewood rearrangement, Pólya-Szegö, and the centered Isoperimetric inequality. Based on these well-known facts, we deduce some straightforward results, e.g., volume comparison and monotonicity of certain functions. 		
				
		\item[Section 3:] Contains the introduction of function spaces and precise definitions of suitable Lorentz-spaces and shows that these spaces are endowed with a natural \emph{Lorentz-type quasinorm}. Moreover, we prove several embedding theorems on these spaces, which are new in the literature.

		\item[Section 4:] This section is devoted to the stability analysis of the Hardy inequality. We begin with classifying all the ``virtual extremals" of the Hardy inequalities on the manifolds. It contains the precise statement of Theorem~\ref{intro-main-th-hardy-cf-stab} and its proof. 
		
		\item[Section 5:] This section is devoted to the study of Sobolev–Lorentz embeddings on Cartan–Hadamard model manifolds. Before presenting and proving the main theorem, we first establish a series of preliminary lemmas and prove Theorem~\ref{avas-best}. 
			
		\item[Section 6:] In this final section, we establish the quantitative equivalence between the Hardy deficit on Cartan–Hadamard model manifolds and its counterpart in Euclidean space. Precise formulations of these results are provided in Theorems~\ref{manifold-eu} and~\ref{eu-manifold}.
		
\end{itemize}

\medskip

\section{Preliminary Geometric settings and Auxiliary results} \label{prem}

In this section, we present an overview of the geometric framework in which our main problems are formulated, together with some straightforward results that will be used repeatedly in the sequel. Several proofs are omitted, as they are standard.

\medskip

\subsection{Riemannian model manifolds and geometric preliminaries}

Let $\m$ be an $N$-dimensional Riemannian manifold ($N \geq 3$), called a \emph{model manifold}, whose metric can be written in the form
\begin{equation*}
g = \mathrm{d}r \otimes \mathrm{d}r + \psi(r)^2 g_{\mathbb{S}^{N-1}},
\end{equation*}
where $r$ denotes the geodesic distance from a fixed point $x_0 \in \m$, referred to as the \emph{pole} of the manifold. The vector $\mathrm{d}r$ represents the radial direction, and $g_{\mathbb{S}^{N-1}}$ is the standard metric on the unit sphere $\mathbb{S}^{N-1}$. The smooth function $\psi : [0,\infty) \to [0,\infty)$ completely determines the geometry of $\m$. 

In this setting, every point $x \in \m \setminus \{x_0\}$ can be expressed in polar coordinates $(r,\theta) \in (0,\infty) \times \mathbb{S}^{N-1}$, where $r$ measures the distance to the pole $x_0$ and $\theta$ specifies the direction of the minimizing geodesic connecting $x$ to $x_0$. For general background on model manifolds, we refer to \cite[Section 3.10]{grig}. Such manifolds arise as a significant subclass of warped product spaces (see, for instance, \cite[Section 1.8]{AMR}).

\medskip 

To ensure that the Riemannian metric determined by the radial function $\psi$ is smooth at the pole, one imposes the following regularity requirements:
\begin{equation}\label{psi}
\psi \in C^\infty([0,\infty)), \qquad 
\psi(r) > 0 \ \text{for } r > 0, \qquad 
\psi'(0) = 1, \qquad 
\psi^{(2k)}(0) = 0 \quad \forall\, k \in \mathbb{N} \cup \{0\}.
\end{equation}
These conditions are both necessary and sufficient for the manifold $\m$ to possess a smooth Riemannian structure at the pole. Although in many analytical contexts, the properties $\psi(0)=0$ and $\psi'(0)=1$ alone guarantee that the metric behaves correctly to first order near the pole, the complete set of conditions in \eqref{psi} is required to prevent any loss of regularity.

Throughout this work, we choose $\psi$ to be defined on the full interval $[0,\infty)$, which ensures that the corresponding model manifold is complete and noncompact—the geometric framework in which our analysis will take place. By the Cartan--Hadamard theorem, any such manifold is topologically equivalent to Euclidean space $\mathbb{R}^N$; more precisely, the exponential map centered at any point is a global diffeomorphism.

Classical examples of noncompact Riemannian model manifolds include hyperbolic space $\mathbb{H}^N$ and Euclidean space $\mathbb{R}^N$, corresponding to the model functions $\psi(r) = \sinh r$ and $\psi(r) = r$, respectively.

In terms of the local coordinate system $\{x^i\}_{i=1}^{N},$ one can write
\begin{align*}
g=\sum g_{ij}{\rm d}x^i\dx^j.
\end{align*}
The Laplace-Beltrami operator $\Delta_g$ concerning the metric $g$ is defined as follows
\begin{align*}
\Delta_g:=\sum \frac{1}{\sqrt{\text{det }(g_{ij})}}\frac{\partial}{\partial x^i}\bigg(\sqrt{\text{det }g_{ij}}\:g^{ij}\frac{\partial}{\partial x^j}\bigg),
\end{align*}
where $(g^{ij})=(g_{ij})^{-1}$. Also, we denote $\nabla_g$ as the Riemannian gradient, and for functions $u$ and $v$, we have\begin{align*}
\langle \nabla_g u , \nabla_g v \rangle_g= \sum g^{ij}\frac{\partial u}{\partial x^i}\frac{\partial v}{\partial x^j}.
\end{align*}
For simplicity, we shall use the notation $|\nabla_g f| = \sqrt{g(\nabla_g f,\nabla_g f)}$. Due to our geometry function $\psi$, we can write these operators more explicitly. Namely, for $x\in\m$, we can write $x=(r, \theta)=(r, \theta_{1},\ldots, \theta_{N-1})\in(0,\infty)\times\sn$ and the Riemannian Laplacian of a scalar function $f$ on $\m$ is given by 
\begin{equation*}
	\lapg u (r, \theta)  =
	\frac{1}{(\psi(r))^2} \frac{\partial}{\partial r} \left[ ((\psi (r)))^{N-1} \frac{\partial u}{\partial r}(r, \theta) \right] \\
	+ \frac{1}{(\psi(r))^2 } \Delta_{\mathbb{S}^{N-1}} u(r, \theta),
\end{equation*}
where $\Delta_{\mathbb{S}^{N-1}}$ is the Riemannian Laplacian on the unit sphere $\mathbb{S}^{N-1}$. Also, let us recall that the Gradient in terms of the polar coordinate decomposition is given by
\begin{equation*}
	\gradg u(r, \theta)=\bigg(\frac{\partial u}{\partial r}(r, \theta), \frac{1}{\psi(r)}\nabla_{\sn}u(r, \theta)\bigg),
\end{equation*}
where $\nabla_{\sn}$ denotes the Gradient on the unit sphere $\sn$.

For any $\delta>0$, we denote by $B_\delta(o):=\{x\in \m\,:\, \text{dist}(x,x_0)<\delta \}$ the geodesic ball in $\m$ with center at a fixed pole $x_0$ and radius $\delta$. For any function $u\in L^1(\m)$, the polar coordinate decomposition can be written as follows
\begin{align*}
\int_{\m}u(x)\dvg=\int_{\sn}\int_{0}^{\infty}u(r,\theta)\:(\psi(r))^{N-1}\dr\dsn,
\end{align*}
where $\dsn$ denotes the measure on the unit sphere.


\subsection{Further geometric assumptions on \texorpdfstring{$\m$}{M} and some intermediate results}

Throughout this article, we additionally assume that $\psi$ is a \emph{convex} function. Under this assumption, the manifold $\m$ becomes a Cartan--Hadamard manifold, i.e., a complete, simply connected Riemannian manifold with nonpositive sectional curvature. For further details, we refer the reader to 
\cite{GW}. In particular, $\psi$ satisfies

\begin{align}\label{convex-assump}
    \psi^{\prime \prime}(r)\geq 0 \quad \text{in } (0,+\infty),\,
\end{align}
which in turn implies 
\begin{align}\label{ins-assump}
    \psi'(r)\geq 1 \quad \text{in } (0,+\infty).\,
\end{align}

It is known that there exists an orthonormal frame $\{F_{j}\}_{j = 1, \ldots, N}$ on $\m$, where $F_{N}$ corresponds to the radial coordinate and $F_{1}, \ldots, F_{N-1}$ correspond to the spherical coordinates, such that $F_{i} \wedge F_{j}$ diagonalize the curvature operator $\mathcal{R}$:  
\[
	\mathcal{R}(F_{i} \wedge F_{N}) = - \frac{\psi^{\prime \prime}}{\psi} \, F_{i} \wedge F_{N}, \quad i < N,
\]
\[
	\mathcal{R}(F_{i} \wedge F_{j}) = - \frac{(\psi^{\prime})^2 - 1}{\psi^2} \, F_{i} \wedge F_{j}, \quad i,j < N.
\]
The quantities
\begin{equation*}
		K_{rad} := - \frac{\psi^{\prime \prime}}{\psi} \quad \text{and} \quad K_{\perp} := - \frac{(\psi^{\prime})^2 - 1}{\psi^{2}}
\end{equation*}
coincide, respectively, with the sectional curvatures of planes containing the radial direction and of planes orthogonal to it.


Assumption \eqref{convex-assump} led us to an immediate straightforward lemma: 

\begin{lemma}
    Suppose $\psi$ satisfies \eqref{psi} and \eqref{convex-assump}. Then the following holds: 
  
\begin{align}\label{sturm}
 \psi(r) \geq r, \quad \mbox{and} \quad  \frac{\psi'(r)}{\psi(r)} \geq \frac{1}{r} \quad \text{for all } r \in (0,+\infty). 
 \end{align}
  \begin{align*}
  \frac{\psi(r_1)}{\psi(r_2)} \geq \frac{r_1}{r_2} \quad \text{for all } \quad r_1>r_2>0.
\end{align*}
\end{lemma}

\medskip

\begin{lemma}\label{aux-fnl}
 Consider a family
 \begin{align*}
     \mathcal{F}:=\{\psi:[0,\infty)\to \mathbb{R} \, :\, \psi \text{ satisfies } \eqref{psi}, \text{ and } \eqref{convex-assump}\}.
 \end{align*}
Let $\psi_0\in\mathcal{F}$. Then $\frac{\psi}{\psi_0}$ is an increasing function for every $\psi\in \mathcal{F}$ if and only if $\psi_0(r)\equiv r$.
\end{lemma}


\begin{proof}
Assume $\frac{\psi}{\psi_0}$ is increasing which in turn implies 
$ \frac{\psi'(r)}{\psi(r)}\geq \frac{\psi_0'(r)}{\psi_0(r)} 
$
for all $r>0.$ In particular,  $\psi(r)=r$ is a member of $\mathcal{F}$ and using \eqref{sturm} for $\psi_0$, it follows that 

\[
\frac{\psi_0'(r)}{\psi_0(r)}= \frac{1}{r} \qquad \text{for all}\quad\, r>0.
\]
Upon integrating above and applying $\psi_0'(0)=1$, we conclude $\psi_0(r)=r$. Conversely, we notice $\big(\frac{\psi(r)}{r}\big)'=\frac{r\psi'(r)-\psi(r)}{r^2}$. Thereby, using the above Lemma, we conclude the proof.
\end{proof}


\begin{remark}
{\rm 
Also, we know that for this manifold $\text{Cut}(x_0)=\emptyset$ (see reference \cite{grig}), where $x_0$ is any fixed point of $\m$ and so, the inverse image of the set $\{0\}$ under the distance function is only the pole $\{x_0\}$, that is $r^{-1}(0)=\{x_0\}$.}
\end{remark}

\medskip

\begin{lemma}\label{complete-lem}
Let $N\geq 3$. Suppose $\psi$ satisfies \eqref{psi} and  \eqref{sturm}. Then the following holds:  \begin{align}\label{complete-lem-eq}
    \sup_{r>0} \frac{\left(r^3\psi(r)^{N-1}\right)^{{\frac{N-2}{N+2}}}}{\int_0^r \frac{\psi(t)^{N-1}}{t^2}\dt}<+\infty.
\end{align}
\end{lemma}
\begin{proof}
First, we notice that due to \eqref{psi}, there holds $\psi(r)=r+O(r^3)$ for $r\to 0^+$. Using this fact, we deduce
 \begin{equation}\label{complete-lem-eq-1}
         \lim_{ r\to 0^+}\frac{\left(r^3\psi(r)^{N-1}\right)^{{\frac{N-2}{N+2}}}}{\int_0^r \frac{\psi(t)^{N-1}}{t^2}\dt}=(N-2).
\end{equation}

Now using $N\geq 3$, and $\psi(r)\geq r$, we deduce
\begin{align*}
    \lim_{r\to +\infty} \left(r^3\psi(r)^{N-1}\right)^{\frac{N-2}{N+2}} =+\infty \qquad \text{ and } \qquad \lim_{r\to +\infty}\int_0^r \frac{\psi(t)^{N-1}}{t^2}\dt=+\infty.
\end{align*}
Hence, applying L'Hospital's rule, we have
\begin{align*}
     \lim_{ r\to +\infty}\frac{\left(r^3\psi(r)^{N-1}\right)^{{\frac{N-2}{N+2}}}}{\int_0^r \frac{\psi(t)^{N-1}}{t^2}\dt}= \lim_{ r\to +\infty} \frac{\left[\left(r^3\psi(r)^{N-1}\right)^{{\frac{N-2}{N+2}}}\right]'}{\left[\int_0^r \frac{\psi(t)^{N-1}}{t^2}\dt\right]'}.
\end{align*}

Now, after computation, we have
\begin{align*}
  \frac{\left[\left(r^3\psi(r)^{N-1}\right)^{{\frac{N-2}{N+2}}}\right]'}{\left[\int_0^r \frac{\psi(t)^{N-1}}{t^2}\dt\right]'}=3\left(\frac{N-2}{N+2}\right)\left(\frac{r}{\psi(r)}\right)^\alpha+\frac{(N-2)(N-1)}{(N+2)}\left(\frac{r}{\psi(r)}\right)^{\alpha+1}\psi'(r),
\end{align*}
where $\alpha=\frac{4(N-1)}{N+2}$. For large $r$, the first term is bounded due to $\psi(r)\geq r$ and $\alpha>0$. Therefore, it remains to check the behavior of the second factor for large $r$. We denote 
\begin{align*}
E(r):=\frac{\psi'(r)}{f(r)^{\alpha+1}} \qquad \text{ where } \quad   f(r)=\frac{\psi(r)}{r}.  
\end{align*}

If possible assume 
\begin{align*}
    \lim_{r \to +\infty} E(r)=+\infty.
\end{align*}
So there exists a large $R>0$ such that $E(r)\geq 2$ for all $r> R$ and for such $r$ this implies
\begin{align*}
    \psi'(r)\geq 2f(r)^{\alpha+1} \qquad \text{ for all } r\geq R.
\end{align*}
Therefore, for $r\geq R$ and using $\psi(r)\geq r$, we have
\begin{align*}
          f'(r)=\frac{\psi'(r)-f(r)}{r}\geq \frac{2f(r)^{\alpha+1}-f(r)}{r}\geq \frac{2f(r)^{\alpha+1}-f(r)^{\alpha+1}}{r}=\frac{f(r)^{\alpha+1}}{r}.
\end{align*}
Next, integrating from $R$ to an arbitrary $r> R$, we obtain
      \begin{equation*}
          \int_R^r \frac{f'(t)}{f(t)^{\alpha+1}}\dt\geq \int_R^r\frac{1}{t}\dt
      \end{equation*}
and using this, along with $\alpha>0$, we get 
     \begin{equation*}
         \frac{1}{\alpha}f(R)^{-\alpha} \geq \frac{1}{\alpha}\left(f(R)^{-\alpha}-f(r)^{-\alpha}\right) \geq \ln\left(\frac{r}{R}\right).
     \end{equation*}
But this gives a contradiction by taking $r\to+\infty$ as $\ln\left(\frac{r}{R}\right)\to +\infty$. Hence, $ \lim_{r \to +\infty} E(r)<+\infty$ and so 
\begin{align}\label{complete-lem-eq-2}
    \lim_{ r\to +\infty}\frac{\left(r^3\psi(r)^{N-1}\right)^{{\frac{N-2}{N+2}}}}{\int_0^r \frac{\psi(t)^{N-1}}{t^2}\dt}<+\infty.
\end{align}
Therefore, combining \eqref{complete-lem-eq-1} and \eqref{complete-lem-eq-2}, we conclude \eqref{complete-lem-eq}.
\end{proof}

\medskip 

\section{Schwarz symmetrization, Rearrangement invariant spaces, and Embeddings}\label{sec-symm}
This section collects the functional-analytic tools we will use to study quantitative stability for the Hardy inequality on Cartan--Hadamard manifolds. Our treatment highlights two themes: rearrangement and symmetrization adapted to the geometry of the manifold, and  Lorentz-type spaces that incorporate the manifold's volume profile. These ingredients allow us to compare arbitrary functions on the manifold with radially symmetric profiles on an appropriate model space and to transfer sharp Euclidean-type inequalities (and their stability properties) to the geometric setting.

\subsection{Schwarz rearrangements and P\'olya-Szeg\"o inequality.} We begin by recalling the notion of rearrangements on a manifold $\m$ and the basic measure-theoretic facts that underlie symmetrization. For a measurable function 
$u: \m \rightarrow \mathbb{R}$ we denote by 
\[
\mu(t) := |( \{ x \in \m : |u(x)| > t\})|,
\]
is its distribution function, where $|A|:= \mbox{Vol}_g(A)$ and we call $u$ to be \textit{admissible} on $\m$ if 
\begin{equation}\label{admissible}
\mu(t)<+\infty \qquad \text{ for all } \quad t>\text{essinf}\,u.
\end{equation}

Throughout this paper, we will work with admissible functions (see \cite[Eqn~1.1, p.~21]{bnstn} ). Let $u^\star$ is the standard one-dimensional decreasing \emph{Hardy-Littlewood} rearrangement in $\m$, which is defined by 
\begin{align}\label{star}
   u^\star(s)= \inf\{t>0: |\{x\in \m: |u(x)|>t\}|\leq s\} \quad \text{ for all } s\geq 0.
\end{align}

Let us denote the volume of a geodesic ball $B_r(x_0)$ centered at the pole $x_0$ and with radius $r>0$ by the following identity,
\begin{align}\label{vol}
    V(r):=V(B_r(x_0))=\omega_N\int_0^r (\psi(t))^{N-1}\dt,
\end{align}
where $\omega_N$ denotes the  total area of the $(N-1)$-dimensional unit sphere $\sn$. We note that $V'(r)=\omega_N(\psi(r))^{N-1}$, and using \eqref{psi}, we can say that $V'(r)>0$ for all $r\in(0,\infty) $. Therefore, $V$ is strictly increasing. Moreover, it is a smooth and invertible function. Let us denote the inverse of $V$ by $G:(0,\infty)\to (0,\infty)$ that is 
\begin{align}\label{inv-vol}
    G(s):=V^{-1}(s) \qquad \text{ for all } s>0.
\end{align}

Then we can define the \emph{Schwarz rearrangement} of $u$ as the radial function $u^\sharp$ defined for every $x\in \m$ by
\begin{align}\label{sharp}
    u^\sharp(x):=u^\star(V(r)),
\end{align}
where $r=\text{dist}(x,x_0)$.

We say that a noncompact Riemannian model manifold satisfies the \emph{P\'olya–Szeg\"o} inequality if for every $u \in D^{1,2}(\m)$ one has
\begin{align}\label{ps}
\int_{\m}|\gradg u^\sharp|^2\dvg \leq \int_{\m}|\gradg u|^2\dvg.
\end{align}

It is well known that standard manifolds such as $\hn$ and $\rn$ admit the P\'olya–Szeg\"o inequality. However, this property does not necessarily extend to all noncompact Riemannian manifolds. Thanks to the recent work of Muratori and Volzone \cite[Theorem~3.9]{mv}, explicit examples are now available where \eqref{ps} fails. At first sight, this might suggest that the class of models satisfying \eqref{ps} is rather restrictive. Nevertheless, this is not the case. In the same work \cite{mv}, the authors established sufficient conditions ensuring the validity of \eqref{ps}. In particular, they proved that for any model manifold, \eqref{ps} holds whenever $\m$ supports a centered isoperimetric inequality, namely,
\begin{align}\label{isoper}
\text{Per}(B_r(x_0)) \leq \text{Per}(\Omega) \quad \forall \,  \Omega \in \mathcal{B}_b(\m),
\end{align}
where $\text{Vol}_g(B_r(x_0)) = V(r) = \text{Vol}(\Omega)$ and $\mathcal{B}_b(\m)$ denotes the family of all bounded Borel sets in $\m$ (see \cite[Proposition~3.5]{mv}).   Apart from $\rn$, it is worth mentioning that \eqref{isoper} holds for $\hn$ (see \cite{Schmidt}) and for $\mathbb{S}^N$ (see \cite{ps-ball}). Moreover, the quantitative version of the centered isoperimetric inequality on $\hn$ was proved in \cite{ps-hn}.  Whether \eqref{isoper} holds in any other model manifolds is not known so far. In a recent preprint \cite{avas-paper}, the first author of this article establishes several necessary conditions and a sufficient condition of the validity of  \it centered isoperimetric inequality\rm. Since this inequality will serve as one of the fundamental tools in our analysis, whenever we require \eqref{ps}, we shall restrict our attention to Riemannian model manifolds satisfying \eqref{isoper}.
Exploiting \eqref{sturm}, we prove the following:

\begin{lemma}
Let $N \geq 2,$ and assume that $\psi$ satisfies assumptions \eqref{psi} and \eqref{convex-assump}. Then there holds 
\begin{align}\label{ve}
    V(r)\geq \frac{\uv}{N}\, r^N \quad \text{for all } r \in (0,+\infty),
\end{align}
and
\begin{align}\label{ge}
    G(s)\leq \left(\frac{N}{\uv}\right)^{\frac{1}{N}}s^{\frac{1}{N}} \quad \text{for all } s \in (0,+\infty).
\end{align}
\end{lemma}


\subsection{Rearrangement invariant subspaces} 

For $p\in (1,\infty)$ and $q\in[1,\infty]$ we define the \textit{Lorentz-type quasinorm} on $\m$ with $N\geq 1$ by,
\begin{equation}\label{quasinorm}
  \|u\|_{\tilde{L}^{p,q}(\m)}:=
\begin{cases}
    \left(\int_0^{\infty}\left(\left(G(s)\right)^{\frac{N}{p}}u^\star(s)\right)^{q}\frac{\ds}{s}\right)^{\frac{1}{q}}\,\,\,\text{if}\,\,\, q< \infty;\\
    \sup_{s>0} \left(G(s)\right)^{\frac{N}{p}}u^\star(s)\,\,\,\text{if}\,\, q=\infty,
\end{cases}  
\end{equation}
where $G$ and $u^\star$ are defined in \eqref{inv-vol} and \eqref{star}, respectively.

\medskip 

We now establish a sub-multiplicative property for the inverse of the volume function, as stated in the following lemma.

\begin{lemma}\label{doubling}
Let $N \geq 2$ and $\psi$ be a positive and increasing function. Then, for all $s>0$ and $\kappa\geq 1$, we have the following:
\begin{align*}
    G(\kappa \,s) \leq \kappa \, G(s). 
\end{align*}
\end{lemma}
\begin{proof}
Fix $r> 0$. Since $\psi$ is increasing, using the change of variable, $t=\kappa x$ we can write 
\begin{align*}
\int_0^{\kappa\,r}\psi(t)^{N-1}\dt = \kappa\int_0^r \psi(\kappa\,x)^{N-1}\dx \geq \kappa\int_0^r \psi(x)^{N-1}\dx.
\end{align*}
Hence, we deduce $V(\kappa\,r) \geq \kappa\,V(r)$ for any $r>0$. Now, as $G$ and $V$ are inverses of each other, taking $r=G(s)$ and so $s=V(r)$, we obtain $V(\kappa\,r) \geq \kappa\,s$. Using the monotonicity of these functions, we get $G(\kappa\,s) \leq \kappa\, G(s)$, as required.
\end{proof}

\medskip

The following proposition shows that $\tilde{L}^{p,q}(\m)$ is a quasinormed space for every $p \in (1,\infty)$ and $q \in [1,\infty]$. It is worth noting that, in the one-dimensional setting, this quasinorm agrees with the usual Lorentz norm on the real line up to a positive multiplicative constant.

\begin{proposition}
Let $(\m,g)$ be an $N$-dimensional Cartan–Hadamard model manifold with $N\geq 2$. Then $\tilde{L}^{p,q}(\m)$ is a quasinormed space for every $p\in (1,\infty)$ and $q\in [1,\infty]$.
\end{proposition}
\begin{proof}
First, we consider the case $q\in[1,\infty)$. If $u=0$, then $u^{\star} \equiv 0$, which trivially implies that $\|u\|_{\tilde{L}^{p,q}(\m)}=0$. Conversely, suppose that $\|u\|_{\tilde{L}^{p,q}(\m)}=0$. Then, using the positivity of the function $G$, together with the right-continuity and monotone decreasing property of $u^{\star}(s)$, we conclude that $u=0$ almost everywhere. This establishes the definiteness property of the quasi-norm. 

Next, we turn to prove homogeneity. For any $\alpha\in \mathbb{R}$, the property $\|\alpha u\|_{\tilde{L}^{p,q}(\m)}=|\alpha|\|u\|_{\tilde{L}^{p,q}(\m)}$ follows by using a standard fact
\begin{align*}
 (\alpha u)^\star(s)=|\alpha|u^\star(s), \qquad \text{ for all } \alpha\in\mathbb{R}.   
\end{align*}

Finally, we will prove the quasi-triangle inequality for the functions $u$ and $v$. First, using the definition \eqref{star}, for $s>0$, there holds
\begin{align*}
    (u+v)^\star(s)\leq u^\star \left(\frac{s}{2} \right)+v^\star \left(\frac{s}{2}\right).
\end{align*}

Now using the above  in the definition \eqref{quasinorm} and applying the change of variable $s=2t$, we deduce
\begin{equation*}
        \|u+v\|^q_{\tilde{L}^{p,q}(\m)}\leq 2^{q-1}\int_0^{\infty} \left[G(2t)^{\frac{N}{p}}u^\star(t)\right]^q\frac{\dt}{t}+2^{q-1}\int_0^{\infty} \left[G(2t)^{\frac{N}{p}}v^\star(t)\right]^q\frac{\dt}{t}.
\end{equation*}
Then, exploiting Lemma~\ref{doubling} in the above, we compute
\begin{equation*}
        \|u+v\|^q_{\tilde{L}^{p,q}(\m)}\leq 2^{q\left(\frac{N}{p}+1\right)-1}\left(\int_0^{\infty} \left[G(t)^{\frac{N}{p}}u^\star(t)\right]^q\frac{\dt}{t}+\int_0^{\infty} \left[G(t)^{\frac{N}{p}}v^\star(t)\right]^q\frac{\dt}{t}\right).
\end{equation*}
Therefore, applying $(a+b)^{\frac{1}{q}}\leq a^{\frac{1}{q}} +b^{\frac{1}{q}}$ for all $a,\, b\geq 0$ and using the definition of Lorentz type quasinorm \eqref{quasinorm}, we get
\begin{equation*}
    \|u+v\|_{\tilde{L}^{p,q}(\m)}\leq 2^{\left(\frac{N}{p}-\frac{1}{q}\right)+1}\left(\|u\|_{\tilde{L}^{p,q}(\m)}+\|v\|_{\tilde{L}^{p,q}(\m)}\right).
\end{equation*}
This completes the proof that  $\tilde{L}^{p,q}(\m)$ is a quasinormed space. The remaining case $q = \infty$ can be treated 
analogously. 
\end{proof}


\medskip

We now establish the following embedding theorem for this Lorentz-type space:
\begin{align*}
    \tilde{L}^{p, q_1}(\m) \hookrightarrow \tilde{L}^{p, q_2}(\m),
\end{align*}
where $p \in (1,\infty)$ and $q_1, q_2 \in [1,\infty]$ with $q_1 \leq q_2$. 
This embedding reflects the usual monotonicity of the Lorentz scale with respect to the secondary index and will be a key ingredient in the functional-analytic framework developed here. 
\medskip

\begin{proposition}
Let $(\m,g)$ be an $N$-dimensional Cartan–Hadamard model manifold with $N\geq 2$. Suppose $p \in (1,\infty)$ and $q_1, q_2 \in [1,\infty]$ with $q_1 \leq q_2$. Then for any measurable function $u$, there holds
\begin{equation}\label{wehave}
   \|u\|_{\tilde{L}^{p,q_2}(\m)}\leq  \left(\frac{2^{\frac{N}{p}}}{\ln 2}\right)^{\frac{q_2-q_1}{q_1q_2}}\|u\|_{\tilde{L}^{p,q_1}(\m)}.
\end{equation}
\end{proposition}
\begin{proof}
The case $q_1=q_2$ is trivial. Next, we may assume $q_2=\infty$ and $q_1<\infty$. Then for any $s>0$, using the decreasing property of $u^\star$, and increasing property of $G$, we observe that 
\begin{align*}
\int_{\frac{s}{2}}^s \left[G(t)^{\frac{N}{p}}u^\star(t)\right]^{q_1}\frac{\dt}{t}\geq \int_{\frac{s}{2}}^s\left[G\left(\frac{s}{2}\right)^{\frac{N}{p}}u^\star(s)\right]^{q_1}\frac{\dt}{t}=\left[G\left(\frac{s}{2}\right)^{\frac{N}{p}}u^\star(s)\right]^{q_1}\, (\ln 2).
\end{align*}
Now, by rearranging the above and using the submultiplicative Lemma~\ref{doubling} with $\kappa =2$, we get
\begin{align*}
   G\left(s\right)^{\frac{N}{p}}u^\star(s)\leq 2^{\frac{N}{p}}\, G\left(\frac{s}{2}\right)^{\frac{N}{p}}u^\star(s)&\leq \left(\frac{2^{\frac{N}{p}}}{\ln 2}\right)^{\frac{1}{q_1}}\left(\int_{\frac{s}{2}}^s \left[G(t)^{\frac{N}{p}}u^\star(t)\right]^{q_1}
   \frac{\dt}{t}\right)^{\frac{1}{q_1}}\\
   &\leq \left(\frac{2^{\frac{N}{p}}}{\ln 2}\right)^{\frac{1}{q_1}}\|u\|_{\tilde{L}^{p,q_1}(\m)}.
\end{align*}
Taking supremum over $s>0$ on the l.h.s. of the above, we obtain 
\begin{align*}
 \|u\|_{\tilde{L}^{p,\infty}(\m)} \leq \left(\frac{2^{\frac{N}{p}}}{\ln 2}\right)^{\frac{1}{q_1}}\|u\|_{\tilde{L}^{p,q_1}(\m)}.
\end{align*}
Next, we are considering the case $q_2<\infty$. Then, using the last case and rearranging, we have
\begin{align*}
    \|u\|^{q_2}_{\tilde{L}^{p,q_2}(\m)}=\int_{0}^\infty \left[G(t)^{\frac{N}{p}}u^\star(t)\right]^{q_2-q_1+q_1}\frac{\dt}{t}&\leq \left( \sup_{s>0} \left(G(s)\right)^{\frac{N}{p}}u^\star(s)\right)^{q_2-q_1}\left( \int_{0}^\infty \left[G(t)^{\frac{N}{p}}u^\star(t)\right]^{q_1}\frac{\dt}{t}\right)\\& = \|u\|_{\tilde{L}^{p,\infty}(\m)}^{q_2-q_1}\|u\|^{q_1}_{\tilde{L}^{p,q_1}(\m)}\\& \leq  \left(\frac{2^{\frac{N}{p}}}{\ln 2}\right)^{\frac{q_2-q_1}{q_1}}\|u\|^{q_2}_{\tilde{L}^{p,q_1}(\m)}.
\end{align*}
Therefore, we have \eqref{wehave}.
\end{proof}

\medskip
We now recall the standard Lorentz quasinorm.

\begin{definition}
{\rm 
 For $p \in (1,\infty)$ and 
$q \in [1,\infty]$, the \textit{Lorentz quasinorm} is defined by
\begin{equation}\label{org-quasinorm}
  \|u\|_{L^{p,q}(\m)}=
\begin{cases}
   \left(\displaystyle\int_0^{\infty}\!\!\left(s^{\frac{1}{p}}\, 
   u^\star(s)\right)^{q}\frac{\mathrm{d}s}{s}\right)^{\frac{1}{q}},
   & \text{if } q<\infty,\\[1.2em]
   \displaystyle\sup_{s>0}\, s^{\frac{1}{p}}u^\star(s),
   & \text{if } q=\infty.
\end{cases}
\end{equation}

We shall show that the Lorentz-type quasinorm introduced earlier 
(see~\eqref{quasinorm}) defines a space that is, in general, larger than 
the classical Lorentz space associated with \eqref{org-quasinorm}. 
Nevertheless, in the special case $\psi(r)=r$, corresponding to the Euclidean 
setting, the two quasinorms agree up to a positive multiplicative constant. 
Indeed, in Euclidean space, one has
\[
V(r)=\frac{\omega_N}{N}r^N,
\qquad\text{and hence}\qquad
G(s)=\left(\frac{N}{\omega_N}\right)^{\frac{1}{N}}s^{\frac{1}{N}},
\]
which shows that both constructions reduce to the usual Lorentz norm 
on $\mathbb{R}^N$.

}
\end{definition}

\medskip

\begin{proposition}
Let $(\m,g)$ be an $N$-dimensional Cartan–Hadamard model manifold with $N\geq 1$. The usual Lorentz space is continuously embedded in the Lorentz-type space. Precisely, for any $p \in (1,\infty)$ and $q \in [1,\infty]$, we have
\begin{align}\label{emb}
    L^{p,q}(\m) \hookrightarrow \tilde{L}^{p,q}(\m).
\end{align}
\end{proposition}
\begin{proof}
The proof immediately follows by arranging the variable and using \eqref{ge}.
\end{proof}

We now recall the classical definition of a \textit{rearrangement-invariant} space, 
as introduced by Luxemburg (see \cite[Definition~4.1, p.~59]{BS} for details).

\begin{definition}
Let $\Omega \subset \m$ be a measurable set.  
A \emph{rearrangement-invariant space} $X(\Omega)$ is a Banach function space endowed with a norm $\|u\|_{X(\Omega)}$ such that
\begin{align*}
    \|u\|_{X(\Omega)} = \|v\|_{X(\Omega)} \quad \text{whenever } u^\star = v^\star.
\end{align*}
\end{definition}

\medskip 

\begin{remark}
{\rm
Some classical examples of rearrangement-invariant spaces include the standard $L^p$ spaces as well as Lorentz spaces. 
From the definition, it can be verified that the Lorentz-type space $\tilde{L}^{p,q}(\m)$ also belongs to this class of spaces. 
In the following, we establish a more precise structural result: the space $\tilde{L}^{2^\star, \infty}(\m)$ is, in fact, the smallest rearrangement-invariant space that contains the entire family of functions $\mathcal{U}_a,$ defined for $N \geq 3,$

\[
\mathcal{U}_a (x) := a(r(x))^{\frac{2-N}{2}}, \quad x \in \m,
\]
where $a \in \mathbb{R} \setminus \{ 0\}$ is a parameter, and $r(x):= \mathrm{dist}_g(x, x_0)$ is the geodesic distance. 
}
\end{remark}

\begin{proposition}
Let $(\m,g)$ be an $N$-dimensional Cartan–Hadamard model manifold with $N\geq 3$.  Then  $\tilde{L}^{2^\star,\infty}(\m)$ is the smallest rearrangement invariant space containing the family of functions $\mathcal{U}_a,$ for $a \in \mathbb{R} \setminus \{ 0\}.$
\end{proposition}
\begin{proof}
    It is enough to work with $v_1(x)=r^{\frac{2-N}{2}}$, where $r:= r(x)$. Now, by simple observation, for any $t>0$, it follows
    \begin{align*}
        |\{x\in \m: |v_1(x)|>t\}|= V\left(t^{-\frac{2}{N-2}}\right).
    \end{align*}
    Therefore, for any $s>0$, we deduce
    \begin{align}\label{ext-rear}
       v_1^\star(s)&= \inf\left\{t>0\,:\, \, V\left(t^{-\frac{2}{N-2}}\right)\leq s \right\}\nonumber \\&= \inf\left\{t>0\,:\,  t^{-\frac{2}{N-2}}\leq G\left( s \right) \right\}= \left(G\left( s \right)\right)^{\frac{2-N}{2}}.
    \end{align}

This implies, 
\begin{align*}
        \|v_1\|_{\tilde{L}^{2^\star,\infty}(\m)}&=\sup_{s>0}\left(G(s)\right)^{\frac{N-2}{2}}v_1^\star(s)=1.
\end{align*}
Finally, this space is the smallest rearrangement-invariant space, which follows from the same argument as given in \cite[Prop.~2.3]{cf-aihp}.
\end{proof}

\medskip 

\begin{remark}\label{only-euclid}
{\rm 
Moreover, we prove that if $\mathcal{U}_a \subset L^{2^\star,\infty}(\m)$ for every $a \in \mathbb{R}$, then the manifold $\m$ must be Euclidean.  Indeed, assuming that $v_1 \in L^{2^\star,\infty}(\m)$, the quantity
\[
\sup_{s>0} \frac{s^{\frac{1}{N}}}{G(s)}
\]
must be finite. Hence, there exists a sequence $\{s_n\}$ with $s_n \to +\infty$ as $n \to +\infty$ such that
\begin{align}\label{exact-ecld}
 s_n^{\frac{1}{N}} \leq C\, G(s_n),   
\end{align}
for some constant $C>0$. In fact, using \eqref{ge} and we have $C\geq \left(\frac{\uv}{N}\right)^{\frac{1}{N}}$. By exploiting \eqref{exact-ecld}, the definition and monotonicity of $G$, we obtain
\[
\int_0^{s_n} \left[\psi(t)^{N-1} - \left(\frac{\uv}{N}\right)\frac{1}{C^N}\,t^{N-1}\right] \,\mathrm{d}t \leq 0
\quad \text{for all } n.
\]
On the other hand, using \eqref{ve} and $\left(\frac{\uv}{N}\right)\frac{1}{C^N}\leq 1$, we also have
\[
\int_0^{s_n} \left[\psi(t)^{N-1} - \left(\frac{\uv}{N}\right)\frac{1}{C^N}\,t^{N-1}\right] \,\mathrm{d}t \geq 0.
\]
Combining these two estimates, we deduce that
\[
\psi(t) = \left(\frac{\uv}{N}\right)\frac{1}{C^N}\, t \quad \text{for every } t \in [0,s_n)
\]
for all $n$. Since $s_n \to +\infty$ as $n \to +\infty$, the above identity extends to all $t \ge 0$. Finally, using the initial condition $\psi'(0) = 1$, we conclude that $C= \left(\frac{\uv}{N}\right)^{\frac{1}{N}}$. Therefore, $\m$ coincides with the Euclidean space precisely when the extremizers belong to the Lorentz space introduced in \eqref{org-quasinorm}. The converse follows immediately from \cite[Prop.~2.3]{cf-aihp}.

}
\end{remark}

\medskip 

In the next remark, we exhibit an example on a Cartan--Hadamard model manifold whose model function grows exponentially at infinity, showing that the Lorentz-type space introduced in \eqref{quasinorm} is strictly larger than the classical Lorentz space defined in \eqref{org-quasinorm}.

\begin{remark}
{\rm
From Remark~\ref{only-euclid}, it follows immediately that $v_1 \notin L^{2^\star,\infty}(\m)$ whenever $\m$ is non-Euclidean. Moreover, it is well known that any one-dimensional nondecreasing and right-continuous function arises as the Hardy--Littlewood rearrangement of some admissible function on a general manifold (see \cite[Chapter~1]{bnstn}). 

Fix $p \in (1,\infty)$ and $q \in [1,\infty)$. Consider an admissible function $u$ whose decreasing rearrangement is given by
\[
u^\star (s)=
\begin{cases}
1, & 0<s<e,\\[4pt]
s^{-\frac{1}{p}}(\ln s)^{-\frac{1}{q}}, & s \geq e,
\end{cases}
\]
on a manifold $(\m,g)$ whose inverse volume function \eqref{inv-vol} satisfies the asymptotic condition
\begin{equation}\label{proper-con}
   G(s)^N = O\!\left(\frac{s}{(\ln s)^\alpha}\right)
   \qquad \text{as } s\to +\infty,
\end{equation}
for some $\alpha>0$. It is clear that \eqref{proper-con} holds in the hyperbolic space $\hn$. More generally, it is valid for any manifold satisfying $e^{\beta r}=O(\psi(r))$ as $r\to+\infty$ for some $\beta>0$.

Under condition \eqref{proper-con}, the function $u$ constructed above belongs to $\tilde{L}^{p,q}(\m)$ but not to $L^{p,q}(\m)$. Indeed,
\[
\|u\|_{L^{p,q}(\m)}^q 
\geq \int_{e}^{\infty} \frac{1}{\ln s}\,\frac{\mathrm{d}s}{s}
= +\infty.
\]
On the other hand, applying \eqref{ge} and \eqref{proper-con}, for large $R>e$, $\frac{p}{q}>0$, and $\alpha>0$, we estimate
\begin{align*}
  \int_0^\infty \left(G(s)^{\frac{N}{p}}\,u^\star (s)\right)^q\frac{\mathrm{d}s}{s} 
  &= \int_0^e\left(G(s)^{\frac{N}{p}}u^\star(s)\right)^q\frac{\mathrm{d}s}{s}
   +\int_e^R \left(G(s)^{\frac{N}{p}}u^\star(s)\right)^q\frac{\mathrm{d}s}{s}+\int_R^\infty\left(G(s)^{\frac{N}{p}}u^\star(s)\right)^q\frac{\mathrm{d}s}{s} \\
  &\leq \frac{N}{\omega_N}\int_0^e s^{\frac{q}{p}-1}\,\mathrm{d}s 
   + \int_e^R \left(G(s)^{\frac{N}{p}}u^\star(s)\right)^q\frac{\mathrm{d}s}{s}
   + C\int_R^\infty(\ln s)^{-1-\frac{\alpha q}{p}}\frac{\mathrm{d}s}{s} \\
  &< +\infty.
\end{align*}
Thus, $u \in \tilde{L}^{p,q}(\m)\setminus L^{p,q}(\m)$, confirming that the space $\tilde{L}^{p,q}(\m)$ is strictly larger than the classical Lorentz space $L^{p,q}(\m)$ in the presence of exponential geometric growth.

}
\end{remark}

\medskip

\begin{remark}
{\rm 
Following the arguments in \cite[Prop.~4.2, p.~217]{BS}, we know that $L^{2^\star}(\m)$ is continuously embedded into its weak counterpart $L^{2^\star,\infty}(\m)$. Therefore, invoking the Sobolev embedding theorem together with \eqref{emb}, we conclude that for any $N$-dimensional Cartan--Hadamard model manifold $(\m,g)$ with $N\geq 3$, the following chain of continuous embeddings holds:
\begin{equation}\label{fnl-emb}
  D^{1,2}(\m)\hookrightarrow L^{2^\star}(\m)\hookrightarrow 
  L^{2^\star,\infty}(\m)\hookrightarrow \tilde{L}^{2^\star,\infty}(\m).
\end{equation}

}
\end{remark}

\medskip

\section{Stability of Hardy Inequality: Classification of Extremals and its stability Via Lorentz-type Embeddings}\label{sec-3}

This section focuses on the quantitative stability of the Hardy inequality on Cartan–Hadamard manifolds. Beyond proving the inequality, we quantify how the deficit controls the distance of a function from the family of extremals. The geometric features of Cartan–Hadamard manifolds—most notably curvature effects on volume growth and radial extremizers—pose substantial analytical difficulties and necessitate tools beyond the Euclidean framework. To overcome these challenges, we rely on the Lorentz-type spaces adapted to model manifolds introduced in the previous section, which provide the functional setting underlying our stability results.

\medskip 
Recall that \(D^{1,2}(\m)\) consists of weakly differentiable functions on \(\m\) whose first-order weak derivatives belong to \(L^2(\m)\); moreover, \(D^{1,2}(\m)\) lies in the admissible class \eqref{admissible}. Cartan--Hadamard manifolds support a Hardy-type inequality on this space. In particular, Carron~\cite{Carron} proved that
\begin{align}\label{hardy}
    \int_{\m} |\nabla_g u|^2 \, \dvg 
    \;\geq\; \frac{(N-2)^2}{4} 
    \int_{\m} \frac{|u|^2}{r^2} \, \dvg,
    \qquad \forall\, u \in D^{1,2}(\m),
\end{align}
where \(r\) denotes the geodesic distance from a fixed pole. The constant \(\frac{(N-2)^2}{4}\) is sharp but never attained in \(D^{1,2}(\m)\), leaving room for refinements via remainder terms or suitable subspaces.

The Hardy inequality originated in Hardy’s study of the discrete \emph{Hilbert inequality} \cite{hardy} and has been extensively developed throughout the $20^{\text{th}}$ century; see \cite{kuf} for a historical account. The first multidimensional Euclidean Hardy inequality was established by Leray \cite{leray}. Important refinements include the addition of nonnegative remainder terms by Brezis–Marcus \cite{bm} and Vázquez–Zuazua \cite{vaz}, and the optimal weighted improvement by Adimurthi–Chaudhuri–Ramaswamy \cite{ANR}. Ghoussoub–Moradifam \cite{GM} introduced a unifying approach via Bessel pairs, while Lam–Lu–Zhang \cite{lu-1} studied geometric Hardy inequalities using general distance functions. Further developments can be found in the preface of Ruzhansky–Suragan \cite{rs-book}.

Extensions of the Hardy inequality to Riemannian manifolds have been widely studied since the pioneering work of \cite{Carron}. Optimal weighted forms and sharp constants were obtained in \cite{KO1, KO, YSK, komber}, while general methods for optimal Hardy weights on noncompact manifolds and spectral approaches were developed in \cite{pinch} and \cite{AK}. Sharp inequalities on Cartan–Hadamard manifolds and alternative Riccati-pair techniques appear in \cite{vhn} and \cite{k-maan}. Hardy-type weights also yield improvements of Poincaré inequalities and Sobolev subspaces in \cite{EDG, BGGP, bgr-cvpde, EGR, GR-23}. Comprehensive treatments can be found in \cite{bar, barbatis, DS, kris-3, kris-1, kris-2, th} and the survey \cite{sandeep-sir}.

The principal goal of this work is to investigate the quantitative form associated with \eqref{hardy} on Cartan--Hadamard model manifolds. In recent years, the problem of establishing stability estimates has attracted considerable attention in the analysis of functional and geometric inequalities. Notable results have been obtained for the Caffarelli--Kohn--Nirenberg inequalities \cite{cfll}, the Hardy--Littlewood--Sobolev inequality \cite{Car17, lu-2,clt-aim-25, clt-maan-26,CLT-arxiv}, the Heisenberg--Pauli--Weyl uncertainty principle \cite{MV21, DN24, dgll, RSY}, the Sobolev inequality \cite{be, f-duke, fn, cf-jems}, and the isoperimetric inequality \cite{maggi-ann}, as well as the Heisenberg uncertainty principle on hyperbolic space \cite{dgll} and the Hardy inequality on Cartan--Hadamard manifolds via scaling of extremizers \cite{RSY}, among others.

Moreover, the sharpness of the constant in \eqref{hardy} can be proved using the family of radial functions
\begin{align*}
    \phi_{\epsilon}(r)=
  \begin{cases}
      r^{\frac{N-2}{2}+\epsilon}\,\,\,\,\text{if}\,\,r\in [0,1]\\
      r^{\frac{2-N}{2}-\epsilon}\,\,\,\,\text{if}\,\,r>1,
  \end{cases}
\end{align*}
where $\epsilon>0$. See, for instance,  \cite[Theorem~2.2]{KO}, where the argument used for hyperbolic space goes exactly similarly for a general Cartan--Hadamard model manifold.

In fact, one can derive an exact remainder term for \eqref{hardy} in the form of a sharp identity (see also \cite{flg}). This identity is obtained by exploiting the framework of Bessel pairs---specifically, by considering the pair 
\[
\Big(r^{\,N-1},\; r^{\,N-1}\frac{(N-2)^2}{4}\frac{1}{r^2}\Big)
\]
on the interval \((0,\infty)\), whose associated positive solution is given by \(r^{\frac{2-N}{2}}\); see also \cite[Theorem~1.1]{FLLM}. As a result, for every function \(u \in C_c^{\infty}(\m\setminus\{x_0\})\), one can verify the following Hardy identity:
\begin{multline}\label{eq-h-model}
\int_{\m} |\nabla_g u|^2 \, \dvg 
- 
\frac{(N-2)^2}{4} \int_{\m} \frac{|u|^2}{r^2} \, \dvg
-
\frac{(N-1)(N-2)}{2} \int_{\m} 
\frac{(r\psi' - \psi)}{\psi}\,\frac{|u|^2}{r^2} \, \dvg \\
= \int_{\m} r^{\,2-N}
\left|
\nabla_g\left(r^{-\frac{2-N}{2}}u\right)
\right|^2 \dvg.
\end{multline}

The above identity quantifies precisely how much \eqref{hardy} fails to be attained in the space \(D^{1,2}(\m)\). In particular, since the right-hand side of \eqref{eq-h-model} is always nonnegative and vanishes only when 
\[
\nabla_g\left(r^{-\frac{2-N}{2}}u\right) \equiv 0,
\]
it follows a fortiori that equality in \eqref{hardy} cannot be achieved by any nontrivial function. Therefore, it is natural to investigate the \emph{virtual extremizers}  of the Hardy identity \eqref{eq-h-model}, namely those functions for which the left-hand side vanishes formally. Motivated by this observation, we next introduce the corresponding family of virtual extremizers and analyze their role in the stability problem for \eqref{hardy}.

\medskip

\begin{definition}
We say $u$ is a virtual extremizer of \eqref{eq-h-model} if $u\in H^1_{loc}(\m\setminus\{x_0\})$ and 
    \begin{align}\label{formal-ex-h}
    \int_{\m}r^{2-N} \left\langle \gradg \varphi, \nabla_{g}\left(r^{-\frac{2-N}{2}} u\right)\right\rangle_g\dvg=0 \qquad \text{ for all } \varphi \in C_c^\infty(\m\setminus\{x_0\}).
    \end{align}
\end{definition}

\medskip

\begin{proposition}
Suppose $u$ is a virtual extremizer of \eqref{eq-h-model} , then 
\begin{align*}
    u(x)=C r^{\frac{2-N}{2}} \qquad \text{for a.e. } x\in \m\setminus\{x_0\},
\end{align*}
where $r=\mathrm{dist}_g(x,x_0)$ and $C$ is some constant.
\end{proposition}
 \begin{proof}
 Set
 \[
 h:=r^{-\frac{2-N}{2}}u,\qquad \text{ and } \qquad w:=r^{2-N},
 \]
 so that the weak extremality condition \eqref{formal-ex-h} reads
 \begin{align}\label{ex-1h}
  \int_{\m} w \,\langle \nabla_g \varphi,\nabla_g h\rangle_g \dvg=0
 \qquad\text{for all }\varphi\in C_c^\infty(\m\setminus\{x_0\}).   
 \end{align}

 Fix an arbitrary connected compact set \(K\subset \m\setminus\{x_0\}\). Choose an open set
 \(\Omega\) with \(K\subset\Omega\Subset \m\setminus\{x_0\}\) and a cutoff
 \(\eta\in C_c^\infty(\Omega)\) satisfying \(0\le\eta\le1\) and \(\eta\equiv1\)
 in a neighborhood of \(K\).

 Since \(u\in H^1_{\mathrm{loc}}(\m\setminus\{x_0\})\), the same holds for $h$. Hence, there exists a sequence \(\{h_j\}\subset C_c^\infty(\Omega)\) with \(h_j\to h\) in \(H^1(\Omega)\).
 For each \(j\) put \(\varphi_j:=h_j\eta^2\in C_c^\infty(\Omega)\) and insert
 \(\varphi_j\) into the weak identity \eqref{ex-1h} to obtain
 \[
 0=\int_{\Omega} w\,\langle\nabla_g(h_j\eta^2),\nabla_g h\rangle_g\dvg
 =\int_{\Omega} w\,\eta^2\langle\nabla_g h_j,\nabla_g h\rangle_g\dvg
 +2\int_{\Omega} w\,\eta\,h_j\langle\nabla_g\eta,\nabla_g h\rangle_g\dvg.
 \]
 Letting \(j\to\infty\) and using \(h_j\to h\) and \(\nabla h_j\to\nabla h\) in
 \(L^2(\Omega)\) yields
 \[
 0=\int_{\Omega} w\,\eta^2|\nabla_g h|^2\dvg + 2\int_{\Omega} w\,\eta\,h\langle\nabla_g\eta,\nabla_g h\rangle_g\dvg.
 \]
 By Cauchy--Schwarz and Young's inequality, we have,
 \begin{align}\label{ex-2h}
  \int_{\Omega} w\,\eta^2|\nabla_g h|^2\dvg
 \le 4\int_{\Omega} w\,h^2|\nabla_g\eta|^2\dvg.
 \end{align}

 Now choose a sequence of cutoffs \((\eta_m)\subset C_c^\infty(\Omega)\)
 satisfying \(\eta_m\equiv1\) on \(K\), \(0\le\eta_m\le1\), \(\operatorname{supp}\eta_m\subset\Omega\),
 and \(\|\nabla_g\eta_m\|_{L^\infty(\Omega)}\to0\) as \(m\to\infty\). Applying \eqref{ex-2h} with \(\eta=\eta_m\) and letting
 \(m\to\infty\) gives
 \[
 \int_{K} w\,|\nabla_g h|^2\dvg = 0.
 \]
 Hence \(\nabla_g h=0\) a.e.\ on \(K\); in particular, \(h\) is a.e.\ constant on
 each connected component of \(K\). Consequently, for some $C_K$ , depending on the connected compact set
 \begin{align}\label{ex-3h}
   u(x)=C_K \, r^{\frac{2-N}{2}}\qquad\text{for a.e.\ }x\in K.  
 \end{align}

 Because dimension $N\geq 2$, and as $\m$ is topologically equivalent to the Euclidean space, so $\m\setminus \{x_0\}$ is still connected. Now consider a compact exhaustion $K_n:=\overline{B_n(x_0)}\setminus B_{\frac{1}{n}}(x_0)$ such that $\cup_{n\geq 1}K_n=\m\setminus \{x_0\}$. Now the discussion \eqref{ex-3h} holds for any arbitrary compact subset of \(\m\setminus\{x_0\}\), we conclude that \eqref{ex-3h} holds for each $K_n$. Moreover, due to exhaustion, we must have
\begin{align*}
   u(x)=C_{K_1} \, r^{\frac{2-N}{2}}\qquad\text{for a.e.\ }x\in \m\setminus\{x_0\} 
 \end{align*}
 with some constant $C_{K_1}$, which completes the proof.
 \end{proof}

\medskip

\subsection{Hardy inequality via Schwartz symmetrization}
To proceed, we first rewrite the inequality \eqref{hardy} in a one-dimensional form by restricting to radially decreasing functions and employing their decreasing rearrangement.

\begin{proposition}\label{1D hardy-prop}
   Let $(\m,g)$ be an $N$-dimensional Cartan–Hadamard model manifold with $N\geq 3$. Let $u\in D^{1,2}(\m)$ be a radially decreasing function and take $s=V(r)$, where $V$ is defined in \eqref{vol} and $r=\mathrm{dist}_g(x,x_0)$. Then \eqref{hardy} reads as
   \begin{equation}\label{1D hardy}
       \left(\frac{N-2}{2}\right)^2\int_0^{\infty}\frac{u^\star(s)^2}{G(s)^2}\ds\leq \uv^2\int_0^{\infty} |{u^\star}'(s)|^2\left(\psi \left(G(s)\right)\right)^{2(N-1)}\ds,
   \end{equation}
   where $u^\star$ is one-dimensional decreasing Hardy-Littlewood rearrangement is defined in \eqref{star} and prime denotes the derivative w.r.t. the $s$ variable.
\end{proposition}
\begin{proof}
As $u$ is a radial decreasing function, we can write $u(x)=u^\sharp(x)=u^\star(V(r))=u^\star(s)$, where $u^\sharp$ is defined in \eqref{sharp}. Also, recall $G$ is the inverse of the volume function $V$ and $\ds=\omega_N (\psi(r))^{N-1}\dr$. Then, using polar coordinates and substitution $r=G(s)$, we compute
\begin{align*}
    \int_{\m} \frac{u^2}{r^2}\dvg=\int_0^{\infty}\frac{u^\star(s)^2}{G(s)^2}\ds,
\end{align*}
and
\begin{align*}
\int_{\m} |\gradg u|^2\dvg=\uv\int_0^{\infty} \left|\frac{\partial u}{\partial r}(r)\right|^2(\psi(r))^{N-1}\dr= \uv^2\int_0^{\infty} |{u^\star} '(s)|^2\left(\psi\left(G(s)\right)\right)^{2(N-1)}\ds.
\end{align*}
Finally, using these two identities in \eqref{hardy}, we obtain \eqref{1D hardy}.
\end{proof}

\begin{remark}
{\rm
We consider the family of functions $\mathcal{U}_a$ given by 
\begin{equation}\label{extremizer}
    \mathcal{U}_a(x):= \mathcal{U}_a(r) = ar^{\frac{2-N}{2}},\qquad a\in \mathbb{R}\setminus\{0\}.
\end{equation}
Notably, these functions are not in any $L^p(\m)$. In fact, \eqref{hardy} has no true non-trivial extremizer, which indicates the possibility of improvement of the inequality.  
}
\end{remark}

\medskip

Our ultimate goal is to refine the Hardy inequality \eqref{hardy} by quantifying how far a given function \(u\) is from the virtual extremal family \(\{\mathcal{U}_a\}\). In this direction, we aim to control the \emph{Hardy deficit}
\begin{align}\label{energy-hardy}
    \mathcal{E}(u)
    :=\int_{\m}|\nabla_g u|^2 \, \dvg 
    - \left(\frac{N-2}{2}\right)^2
    \int_{\m}\frac{|u|^2}{r^2}\,\dvg,
\end{align}
from below by a suitable \emph{remainder term}, explicitly expressed in terms of a distance between \(u\) and the manifold of extremizers. Establishing such a bound constitutes the essence of a \emph{stability} result for the Hardy inequality.

\medskip
 For this stability analysis, we introduce a transformation that allows us to relate the geometry of \(\m\) to that of the Euclidean space \(\mathbb{R}^N\). Let \(\tilde{r} > 0\) be arbitrary and recall that \(V(\tilde{r})\) denotes the Riemannian volume of the geodesic ball \(B_{\tilde{r}}(x_0)\), as defined in \eqref{vol}. We now define a Euclidean radius \(\tilde{\varrho} > 0\) such that the Euclidean ball \(B_{\tilde{\varrho}}(0)\subset \mathbb{R}^N\) has the same volume:
\begin{align*}
    \int_0^{\tilde{r}} \psi(t)^{\,N-1} \, \mathrm{d}t
    \;=\;
    \int_0^{\tilde{\varrho}} t^{N-1} \, \mathrm{d}t,
\end{align*}
where \(\psi\) is the model function generating the metric on \(\m\). By inverting this relation, we deduce that for every \(s>0\) there exists a unique \(\varrho>0\) (and vice versa) such that
\begin{align}\label{transform}
    G(s)=\varrho,
\end{align}
where the function \(G\) is the inverse of the model volume map on the Cartan--Hadamard manifold.

\medskip
This volume-preserving correspondence naturally leads us to introduce the Jacobian distortion factor that measures the deviation of the geometry of \(\m\) from the Euclidean one. Specifically, we define the function \(J:[0,\infty)\to [1,\infty)\) by
\begin{align}\label{jacobian}
    J(t)
    :=
    \left(
    \frac{\psi\!\left(\Big(\frac{N}{\omega_{N}}\,t\Big)^{\!\frac{1}{N}}\right)}
    {\left(\frac{N}{\omega_{N}}\,t\right)^{\!\frac{1}{N}}}
    \right)^{\!N-1},
    \qquad \text{with } J(0)=1.
\end{align}
The quantity \(J(t)\) can be interpreted as the ratio between the area elements in the two geometries when expressed in volume coordinates. The fact that \(J(t)\geq 1\) reflects the non-positive curvature assumption underlying Cartan--Hadamard manifolds.

\medskip
This geometric transformation and the associated distortion factor will serve as fundamental tools in our analysis, particularly when passing to rearranged, radially decreasing representatives of functions in \(D^{1,2}(\m)\), and when comparing functional inequalities on \(\m\) with their Euclidean counterparts.

Now we need a simple property of the function $J$ described below.
\begin{proposition}\label{inc-j}
Suppose $\psi$ satisfies \eqref{psi} and \eqref{sturm}. Then the function $J$ defined in \eqref{jacobian} is smooth and increasing on $[0,\infty)$. 
\end{proposition}
\begin{proof}
First using $\psi'(0)=1$ we deduce $\lim_{t\to 0^+} J(t)=1$ and hence the continuity at zero of $J$ follows. Now writing $\tilde{t}=\big(\frac{N}{\uv}t\big)^{\frac{1}{N}}$, we get $J(t)=\big(\frac{\psi(\tilde{t})}{\tilde{t}}\big)^{N-1}$ and using the Lemma~\ref{aux-fnl} increasing property follows.
\end{proof}

Below, we shall prove several lemmas and propositions to establish the stability results. 

\medskip

\textbf{Influence of the Jacobian-type Transformation.} 
The change of variables introduced in \eqref{transform} allows us to reformulate
the $N$-dimensional Hardy inequality \eqref{hardy} through its 
one-dimensional counterpart \eqref{1D hardy}. In particular, this transformation 
converts the original inequality into a \emph{one-dimensional weighted} Hardy 
inequality related to the weighted Euclidean formulation 
(see \cite[Eqn.~2.5]{cf-aihp}). This weighted version in one dimension will serve 
as a fundamental tool in our subsequent analysis.

\medskip

Therefore, it is important to see how all the variables are connected to each other and \eqref{1D hardy} and the extremizer set \eqref{extremizer} takes the form under all these variables. We start with any $u\in D^{1,2}(\m)$ which is radially decreasing function and take $s=V(r)$, where $V$ is defined in \eqref{vol} and $r=\mathrm{dist}_g(x,x_0)$. Now, by using our transformation \eqref{transform},  for any positive $\varrho$, we know there exists a unique positive $s$. Thereby, we define, 
\begin{align}\label{phirho}
    \phi(\varrho):=u^\star(s)  \qquad \text{ for all } \varrho>0.
\end{align}
Now, under the newly defined above function, inequality \eqref{1D hardy}, takes the following form
\begin{equation}\label{1D Etype}
    \left(\frac{1}{2^\star}\right)^2\int_0^{\infty} |\phi(\varrho)|^2\varrho^{-\frac{2}{N}}J(\varrho)\drho\leq \int_0^{\infty}\left(-\phi'(\varrho)\right)^2\varrho^{\frac{2}{N'}}J(\varrho)\drho
\end{equation}
where  $2^\star= \frac{2N}{N-2}$, $N'=\frac{N}{N-1}$, and $J$ is defined in \eqref{jacobian}. In fact, this can be verified by delicate calculation and using $\ds=J(\varrho)\drho$.

Now we will describe how the extremizers take the form under these changes of variables. For any $a \in \mathbb{R} \setminus \{ 0\}$, we have the Hardy virtual extremizer $\mathcal{U}_a(r)=ar^{\frac{2-N}{2}}$. Then, under the Hardy-Littlewood rearrangement (see \eqref{ext-rear}), and using \eqref{transform}, we have 
\begin{align*}
    \mathcal{U}_a^\star(s)=a \left(G\left( s \right)\right)^{\frac{2-N}{2}}= a \left(\frac{N}{\uv}\varrho\right)^{-\frac{1}{2^\star}}.
\end{align*}

In the following result, we perform a stability analysis of the one-dimensional weighted Euclidean Hardy inequality derived in \eqref{1D Etype}. This analysis further yields a quantitative refinement of the one-dimensional inequality \eqref{1D hardy}, which will serve as a key preparatory tool for one of our main results.

\begin{lemma}\label{stab-1d-e}
Let $N\geq 3$ and $N'$ be its H\"older conjugate. Suppose $\psi$ satisfies \eqref{psi}, \eqref{convex-assump}, and \eqref{ins-assump}. Denote 
\begin{equation}\label{distance 1D Etype}
        \delta(\phi)=\inf_{a\geq 0}\frac{\|\phi(\varrho)-a\varrho^{-\frac{1}{2^\star}}\|_{L^{2^\star,\infty}(0,\infty)}}{\left(\int_0^{\infty}\phi(\varrho)^2\varrho^{-\frac{2}{N}}J(\varrho)\drho\right)^{\frac{1}{2}}}.
\end{equation} 
Then
\begin{equation}\label{stability 1D Etype}
        \left(\frac{1}{2^\star}\right)^2\left(\int_0^{\infty} |\phi(\varrho)|^2\varrho^{-\frac{2}{N}}J(\varrho)\drho\right)\left[1+\left(\frac{\delta(\phi)}{2}\right)^4\right]\leq \int_0^{\infty}\left(-\phi'(\varrho)\right)^2\varrho^{\frac{2}{N'}}J(\varrho)\drho,
\end{equation}
holds for every non-increasing locally absolutely continuous function $\phi:(0,\infty)\to[0,\infty)$, making r.h.s. of \eqref{stability 1D Etype} finite and such that $\lim_{\varrho\to \infty}\phi(\varrho)=0$.
\end{lemma}
\begin{proof}
Without loss of generality, we may assume 
\begin{align}\label{c-h-1}
   \int_0^{\infty} |\phi(\varrho)|^2\varrho^{-\frac{2}{N}}J(\varrho)\drho=1.
\end{align}
Now, for the notational economy, we denote
\begin{align}\label{c-h-2}
    \si:= \int_0^{\infty}\left(-\phi'(\varrho)\right)^2\varrho^{\frac{2}{N'}}J(\varrho)\drho.
\end{align}
Therefore, to prove \eqref{stability 1D Etype}, it is enough to establish
\begin{align*}
    \left(\frac{1}{2^\star}\right)^2\left[1+\left(\frac{\delta(\phi)}{2}\right)^4\right]\leq \si.
\end{align*}
Furthermore, denoting
\begin{align}\label{d}
    \sd:=({2^\star}^2\si -1),
\end{align}
we will be through if we establish
 \begin{equation}\label{second}
     \delta(\phi) \leq 2 \,\sd^{\frac{1}{4}}.
\end{equation}

Moreover, it is not harmful to work with $\mathcal{D}>0$. Otherwise, in the case $\mathcal{D}=0$, we can trace back to the equality of the Hardy inequality, and this will give $\delta(\phi)=0$, and in this case, \eqref{stability 1D Etype} follows trivially. Now, we will break the remaining proof into several smaller steps. We will begin with the evaluation of the boundary term near infinity.

\textbf{Step 1.} In this step, we will show  
\begin{align}\label{step-1}
    \lim_{\varrho \to \infty} \left(\phi(\varrho)\right)^2\varrho^{1-\frac{2}{N}}J(\varrho)= 0.
\end{align}
First, for any $\varrho>0$, using $\lim_{\varrho\to \infty}\phi(\varrho)=0$ and  H\"older inequality, we deduce
\begin{align}\label{step-1-1}
    (\phi(\varrho))^2 =\left(\int_\varrho^{\infty}-\phi'(t)\dt\right)^2 \leq \left(\int_\varrho^{\infty}\left(-\phi'(t)\right)^2t^{\frac{2}{N'}}J(t)\dt\right) \left(\int_\varrho^{\infty}\frac{1}{t^{\frac{2}{N'}}J(t)}\dt\right).
\end{align}
Now, for $\varrho>0$, we denote
\begin{align*}
    \si(\varrho):= \int_\varrho^{\infty}\left(-\phi'(\varrho)\right)^2\varrho^{\frac{2}{N'}}J(\varrho)\drho.
\end{align*}
Note that this term $\si(\varrho)$ is always bounded above by the finite term $\si$. Next, using this notation in \eqref{step-1-1}, we deduce 
\begin{align*}
    \left(\phi(\varrho)\right)^2\varrho^{1-\frac{2}{N}}J(\varrho)\leq \si(\varrho) \varrho^{1-\frac{2}{N}}J(\varrho)\int_\varrho^{\infty}\frac{1}{t^{\frac{2}{N'}}J(t)}\dt.
\end{align*}
Now, recall the definition of $J$ defined in \eqref{jacobian}. Then, using the substitution $x= \left(\frac{N}{\uv}\varrho\right)^{\frac{1}{N}}$ and making a change of variable $Nt=\uv y^N$ in the above integration,  we deduce
\begin{align*}
     \varrho^{1-\frac{2}{N}}J(\varrho)\int_\varrho^{\infty}\frac{1}{t^{\frac{2}{N'}}J(t)}\dt = N\frac{(\psi(x))^{N-1}}{x}\int_x^\infty \frac{1}{(\psi(y))^{N-1}}\dy.
\end{align*}
Now, using \eqref{sturm} in the above for the variable $y>x$, we obtain 
\begin{align*}
    \int_x^\infty \frac{1}{(\psi(y))^{N-1}}\dy\leq \frac{x^{N-1}}{(\psi(x))^{N-1}}\int_x^\infty \frac{1}{y^{N-1}}\dy=\frac{1}{(N-2)}\frac{x}{(\psi(x))^{N-1}}.
\end{align*}
Finally, exploiting this in the previous integral estimate, and observing $\lim_{\varrho\to \infty}\si(\varrho)=0$, we deduce \eqref{step-1}.

\textbf{Step 2.} In this step, we will show  
\begin{align}\label{step-2}
    \lim_{\varrho \to 0^+} \left(\phi(\varrho)\right)^2\varrho^{1-\frac{2}{N}}J(\varrho)= 0.
\end{align}
First, we notice that using $J(\varrho)\geq 1$ into \eqref{c-h-1} and \eqref{c-h-2} there holds
\begin{align}\label{step-20}
     \int_0^{\infty} |\phi(\varrho)|^2\varrho^{-\frac{2}{N}}\drho<+\infty, \qquad \text{ and } \quad \int_0^{\infty}\left(-\phi'(\varrho)\right)^2\varrho^{\frac{2}{N'}}\drho<+\infty.
\end{align}
Moreover, mimicking the argument as in Step 1 with $J=1$, we have
\begin{align}\label{step-21}
    \lim_{\varrho \to \infty} \left(\phi(\varrho)\right)^2\varrho^{1-\frac{2}{N}}= 0.
\end{align}

Now for $R>r>0$, the fundamental theorem of calculus yields
\begin{align*}
    \left(\phi(r)\right)^2r^{1-\frac{2}{N}}-\left(\phi(R)\right)^2R^{1-\frac{2}{N}}= -\left(\frac{N-2}{N}\right)\int_r^R (\phi(\varrho))^2\varrho^{-\frac{2}{N}}\drho+2\int_r^R \phi(\varrho)(-\phi'(\varrho))\varrho^{1-\frac{2}{N}}\drho.
\end{align*}
Now, by applying H\"older’s inequality to the last term above and using \eqref{step-20} together with \eqref{step-21}, we obtain that the limit $\displaystyle \lim_{\varrho \to 0^+} \left(\phi(\varrho)\right)^2 \varrho^{1-\frac{2}{N}}$ exists. Now, if possible, assume 
\begin{align*}
    \lim_{\varrho \to 0^+} \left(\phi(\varrho)\right)^2\varrho^{1-\frac{2}{N}}=\ell,
\end{align*}
where $\ell>0$. Then there exists a small $\tilde{r}>0$ such that 
\begin{align*}
 \left(\phi(\varrho)\right)^2\varrho^{-\frac{2}{N}}>\frac{\ell}{2\varrho}\qquad \text{ for all } \quad 0<\varrho<\tilde{r}.
\end{align*}
Now taking integration over $0$ to $\tilde{r}$ in the above, we can see this contradicts the first part of \eqref{step-20}. Therefore, we must have
\begin{align*}
    \lim_{\varrho \to 0^+} \left(\phi(\varrho)\right)^2\varrho^{1-\frac{2}{N}}=0.
\end{align*}
Finally, using this along with $\lim_{\varrho\to 0^+}J(\varrho)=1$, we conclude \eqref{step-2}.

\textbf{Step 3.} In this step, we estimate
\begin{equation}\label{fourth}
    2^\star\int_0^{\infty}\left(-\phi'(\varrho)\right)\phi(\varrho)\varrho^{1-\frac{2}{N}}J(\varrho)\drho\geq 1.
\end{equation}
To establish this, first, we notice
\begin{align*}
        \left(\left(\phi(\varrho)\right)^2\varrho^{1-\frac{2}{N}}J(\varrho)\right)'=2\phi(\varrho)\phi'(\varrho)\varrho^{1-\frac{2}{N}}J(\varrho)+\left(\phi(\varrho)\right)^2\left(\left(1-\frac{2}{N}\right)\varrho^{-\frac{2}{N}}J(\varrho)+\varrho^{1-\frac{2}{N}}J'(\varrho)\right).
\end{align*}
Then, we apply again the fundamental theorem of calculus for the integral in $0<r<R<\infty$. Now letting $r\to 0^+$ and $R\to \infty$, and using the boundary condition derived in Step 1 and Step 2, we prove
\begin{multline}\label{third}
        \int_0^{\infty} |\phi(\varrho)|^2\varrho^{-\frac{2}{N}}J(\varrho)\drho=2^\star\int_0^{\infty}\left(-\phi'(\varrho)\right)\phi(\varrho)\varrho^{1-\frac{2}{N}}J(\varrho)\drho-\frac{N}{N-2}\int_0^{\infty}\phi(\varrho)^2\varrho^{1-\frac{2}{N}}J'(\varrho)\drho.
\end{multline}
Now using Proposition~\ref{inc-j}, we have $J'(\varrho)>0$, hence using this in the above along with \eqref{c-h-1}, we establish \eqref{fourth}.

\textbf{Step 4.} In this part, we will prove
\begin{equation}\label{fifth}
    \int_0^{\infty}\left(A(\varrho)-B(\varrho)\right)^2\drho\leq\mathcal{D},
\end{equation}
where $\mathcal{D}$ is defined in \eqref{d} and for $\varrho>0$, we denote
\begin{align*}
    A(\varrho):= 2^\star\left(-\phi'(\varrho)\right)\varrho^{\frac{1}{N'}}\sqrt{J(\varrho)} \qquad \text{ and } \quad
    B(\varrho):=\phi(\varrho)\varrho^{-\frac{1}{N}}\sqrt{J(\varrho)}.
\end{align*}
Next, rewriting \eqref{third} under these notations and using \eqref{c-h-1} we have
\begin{equation*}
    \int_0^{\infty}A(\varrho)B(\varrho)\drho=1+\frac{N}{N-2}\int_0^{\infty}\phi(\varrho)^2\varrho^{1-\frac{2}{N}}J'(\varrho)\drho.
\end{equation*}
Now, using \eqref{c-h-1}, \eqref{c-h-2}, $J'(\varrho)\geq 0$,  and rearranging we deduce
\begin{equation*}
    \int_0^{\infty}\left(A(\varrho)-B(\varrho)\right)^2\drho=(2^\star)^2 \si-1-2^\star\int_0^{\infty}\phi(\varrho)^2\varrho^{1-\frac{2}{N}}J'(\varrho)\drho,
\end{equation*}
which entails \eqref{fifth}.

\textbf{Step 5.} In this step, we define a key function by
\begin{align}\label{tau}
     \tau(\varrho):=\varrho^{\frac{1}{2^\star}}\phi(\varrho) \qquad \text{ for } \quad \varrho>0.
\end{align}
By direct calculation, we write
\begin{equation*}
    A(\varrho)-B(\varrho)=-2^\star \varrho^{\frac{1}{2}}\tau'(\varrho)\sqrt{J(\varrho)}
\end{equation*}
and using this \eqref{fifth} becomes
\begin{equation}\label{sixth}
\int_0^{\infty}\varrho\left(\tau'(\varrho)\right)^2J(\varrho)\drho\leq\frac{\mathcal{D}}{(2^\star)^2}.
\end{equation}

\textbf{Step 6.} In this final step, we introduce the Marcinkiewicz distance function. For that, we consider the sub-level set by
\begin{equation*}
    \mathcal{A}:=\left\{\varrho>0\,: \tau(\varrho)> \mathcal{D}^{\frac{1}{4}}\right\}.
\end{equation*}

Using $\lim_{\varrho\to 0^+}J(\varrho)=1$, \eqref{step-2}, $J(\varrho)\geq 1$, and \eqref{step-1}, we comment that the function $\tau$ vanishes near zero and infinity. With the help of this observation, we immediately notice that $\mathcal{A}$ is a bounded and open subset of $(0,\infty)$. Hence, it is a countable disjoint union of open, bounded intervals. Let $(a,b)$ be any such connected component of $\mathcal{A}$. On $(a,b)$ we have $\tau(\varrho)>\mathcal{D}^{\frac{1}{4}}$ and exploiting \eqref{tau}, and \eqref{c-h-1}, we evaluate
\begin{align}\label{step61}
    \int_{a}^{b} \frac{J(\varrho)}{\varrho}\drho
    \leq\mathcal{D}^{-\frac{1}{2}}\int_{a}^{b}\frac{\left(\tau(\varrho)\right)^2J(\varrho)}{\varrho}\drho\leq\mathcal{D}^{-\frac{1}{2}}.
\end{align}

From the continuity of the function $\tau$, it follows that at each endpoint of $(a,b)$, the value of $\tau$ is exactly $\mathcal{D}^{\frac{1}{4}}$. Therefore, using this and \eqref{sixth} along with H\"older inequality, for some $\varrho\in (a,b)$, we deduce
\begin{align*}
     |\tau(\varrho)-\mathcal{D}^{\frac{1}{4}}|
    =\left|\int_{a}^\varrho\tau'(t)\dt\right|
     \leq\left(\int_{a}^{b}\varrho\left(\tau'(\varrho)\right)^2J(\varrho)\drho\right)^{\frac{1}{2}}\left(\int_{a}^{b}\frac{1}{\varrho J(\varrho)}\drho\right)^{\frac{1}{2}} \leq \frac{{\mathcal{D}}^{\frac{1}{2}}}{2^\star}\left(\int_{a}^{b}\frac{1}{\varrho J(\varrho)}\drho\right)^{\frac{1}{2}}.
\end{align*}

Now we will estimate the second factor of the above estimate. Using the information $J(\varrho)\geq 1$, and \eqref{step61}, we deduce
\begin{align*}
    \int_{a}^{b}\frac{1}{\varrho J(\varrho)}\drho\leq   \int_{a}^{b}\frac{1}{\varrho }\drho \leq  \int_{a}^{b} \frac{J(\varrho)}{\varrho}\drho \leq \mathcal{D}^{-\frac{1}{2}}.
\end{align*}
Hence, using this in the previous estimate, we have
\begin{align*}
     |\tau(\varrho)-\mathcal{D}^{\frac{1}{4}}|\leq \frac{{\mathcal{D}}^{\frac{1}{4}}}{2^\star}.
\end{align*}
On $\mathcal{A}^c$, by the help of the triangle inequality, we have $|\tau(\varrho)-\mathcal{D}^{\frac{1}{4}}|\leq 2\mathcal{D}^{\frac{1}{4}}$. Combining these two estimates and writing back $\tau$ to $\phi$ from \eqref{tau}, we deduce
\begin{align*}
    \varrho^{\frac{1}{2^\star}}|\phi(\varrho)-\mathcal{D}^{\frac{1}{4}}\varrho^{-\frac{1}{2^\star}}|\leq 2\mathcal{D}^{\frac{1}{4}} \qquad \text{ for all }\varrho>0.
\end{align*}
Therefore, taking the supremum over $\varrho>0$ in the above, and using the definition of Marcinkiewicz space, we deduce
\begin{equation*}
    \|\phi(\varrho)-\mathcal{D}^{\frac{1}{4}}\varrho^{-\frac{1}{2^\star}}\|_{L^{2^\star,\infty}\left(0,\infty\right)}\leq2\mathcal{D}^{\frac{1}{4}}.
\end{equation*}
Hence, 
\begin{equation*}
    \inf_{a\geq0}\|\phi(\varrho)-a\varrho^{-\frac{1}{2^\star}}\|_{L^{2^\star,\infty}\left(0,\infty\right)}\leq2\mathcal{D}^{\frac{1}{4}}
\end{equation*}
and the desired estimate \eqref{second} holds.
\end{proof}

Now we define the Lorentz-type quasinorm on $[0,\infty)$ corresponding to a Cartan--Hadamard model manifold  $(\m,g)$ with dimension $N\geq 3$. Suppose $p\in (1,\infty)$ and $q\in[1,\infty]$. Then for any measurable function $f:[0,\infty)\to \mathbb{R}$, we say it is going to be a member of $\tilde{L}^{p,q}_{\m}(0,\infty)$ if the following is finite
\begin{equation*}
  \|f\|_{\tilde{L}^{p,q}_{\m}(0,\infty)}=
\begin{cases}
    \left(\int_0^{\infty}\left(\left(G(s)\right)^{\frac{N}{p}}f^\star(s)\right)^{q}\frac{\ds}{s}\right)^{\frac{1}{q}}\,\,\,\text{if}\,\,\, q< \infty;\\
    \sup_{s>0} \left(G(s)\right)^{\frac{N}{p}}f^\star(s)\,\,\,\text{if}\,\, q=\infty,
\end{cases}  
\end{equation*}
where $G$  is the inverse of the volume function defined in \eqref{vol} related to $\m$ and and $f^\star$ is defined below
\begin{align*}
   f^\star(s)= \inf\{t>0:  \mathcal{L}^1(\{x\in (0,\infty): |f(x)|>t\})\leq s\} \quad \text{ for all } s\geq 0.
\end{align*}

\medskip

Next, we express the distance function $\delta$, defined in \eqref{distance 1D Etype}, in terms of the variable $s$. Recall that the variables $\varrho$ and $s$ are related through the transformation given in \eqref{transform}. We then proceed to establish the stability analysis of the one-dimensional Hardy inequality for radially decreasing functions.
\begin{proposition}\label{1stabi}
Let $(\m,g)$ be an $N$-dimensional Cartan–Hadamard model manifold with $N \geq 3$, and let $u \in D^{1,2}(\m)$ be a radially decreasing function. Suppose that $u^\star(s)$ and $G(s)$ are defined in \eqref{star} and \eqref{inv-vol}, respectively. Let $\phi(\varrho)$ be the function associated with $u^\star(s)$ through the relation \eqref{phirho}, where the variables $\varrho$ and $s$ are connected via the transformation \eqref{transform}. Denote
\begin{equation}\label{L to Ltilde}
    \nu(u^\star)
    := \inf_{a \geq 0} 
    \frac{
        \left\| u^\star(s) - a \, (G(s))^{\frac{2-N}{2}} \right\|_{\tilde{L}^{2^\star,\infty}_{\m}(0,\infty)}
    }{
        \left( \int_0^{\infty} u^\star(s)^2 \, (G(s))^{-2} \ds \right)^{\frac{1}{2}}
    }.
\end{equation}
Then
\begin{equation}\label{rho to s}
    \delta(\phi) =\left(\frac{\uv}{N}\right)^{\frac{N+2}{2N}} \, \nu(u^\star),
\end{equation}
where $\delta(\phi)$ is defined in \eqref{distance 1D Etype}. Moreover, it holds
\begin{equation}\label{stability 1D Htype}
    \left(\frac{N-2}{2}\right)^2\int_0^{\infty}\frac{u^\star(s)^2}{\left(G(s)\right)^2}\ds\left[1+\left(\frac{\uv}{N}\right)^{\frac{2(N+2)}{N}}\left(\frac{\nu\left(u^\star\right)}{2}\right)^4\right]\leq \uv^2\int_0^{\infty} |{u^\star}'(s)|^2\left(\psi(G(s))\right)^{2(N-1)}\ds.
\end{equation}
\end{proposition}
\begin{proof}
 First, we notice that using \eqref{transform}, \eqref{jacobian} and \eqref{phirho}, we have
 \begin{align}\label{rho to s-1}
    \int_0^{\infty}\phi(\varrho)^2\varrho^{-\frac{2}{N}}J(\varrho)\drho=\left(\frac{\uv}{N}\right)^{-\frac{2}{N}}\int_0^{\infty} u^\star(s)^2 \, (G(s))^{-2} \ds.
 \end{align}
 
 Now, suppose $a\geq 0$, and we define the following functions by
 \begin{align*}
    \tilde{g}(s):=u^\star(s)-a\left(G(s)\right)^{\frac{2-N}{2}} \qquad \text{ and } \quad g(\varrho):= \phi(\varrho)-a\left(\frac{N}{\uv}\right)^{-\frac{1}{2^\star}}\,\varrho^{-\frac{1}{2^\star}}.
 \end{align*}
Since we are working within a rearrangement-invariant framework, it suffices to consider the case where $|\tilde{g}|$ is nonincreasing. Otherwise, one may replace $|\tilde{g}|$ by its decreasing rearrangement without loss of generality. Now we continue the proof by defining the following super-level sets. For $t>0$, we denote
\begin{align*}
        \tilde{E}_t=\left\{s>0:|\tilde{g}(s)|>t \right\} \qquad \text{ and } \quad E_t=\left\{\varrho>0: |g(\varrho)|>t\right\}.
\end{align*}
 Then we see $E_t=\varrho(\tilde{E}_t)=\{\varrho(s):s\in \tilde{E}_t\}$.  Since $|\tilde{g}|$ is a decreasing function, so $\tilde{E}_t=(0,s_t)$, where $\tilde{g}^\star(s_t)=t$ for some $s_t>0$.  Since $\varrho$ is an increasing function w.r.t $s$ variable and using $\varrho(0)=0$, we have $|E_t|=\varrho(s_t)=\frac{\uv}{N}\left(G(s_t)\right)^N$. Therefore, using the equivalent definition of Lorentz space (see \cite[Prop.~1.4.9, p.~53]{grafa}), we have
\begin{align*}
        \|g(\varrho)\|_{L^{2^\star,\infty}(0,\infty)}=\sup_{t>0}t\left(E_t\right)^{\frac{1}{2^\star}}
        =\left(\frac{\uv}{N}\right)^{\frac{1}{2^\star}}\sup_{t>0} \tilde{g}^\star(s_t)\left(G(s_t)\right)^{\frac{N}{2^\star}}.
\end{align*}
Due to $\tilde{g}(s)\neq 0$ on a unbounded set, we have a one-to-one continuous correspondence between $(0,\infty)$ to $(0,\infty)$ between $t$ and $s_t$. Hence, we deduce
\begin{align}\label{rho to s-2}
        \|g(\varrho)\|_{L^{2^\star,\infty}(0,\infty)}
        =\left(\frac{\uv}{N}\right)^{\frac{1}{2^\star}}\sup_{t>0} \tilde{g}^\star(s_t)\left(G(s_t)\right)^{\frac{N}{2^\star}}=\left(\frac{\uv}{N}\right)^{\frac{1}{2^\star}}\|\tilde{g}(s)\|_{\tilde{L}^{2^\star,\infty}_{\m}(0,\infty)}.
\end{align}
Finally, combining \eqref{rho to s-1} and \eqref{rho to s-2}, the relation \eqref{rho to s} follows. Moreover, \eqref{stability 1D Htype}, is an easy consequence of Lemma~\ref{stab-1d-e}.
\end{proof}

\medskip

To carry out the stability analysis for a general function (not necessarily radial), 
we first quantify, in a suitable sense, the distance between the function and its 
Schwarz rearrangement, which is radial by construction. The key idea is that if a 
function is sufficiently close to its rearrangement, then---in light of 
\eqref{stability 1D Htype}---the associated radial representative lies near the set 
of extremizers. Consequently, the original function must also be close to the 
extremizers. Establishing this connection forms the main objective of the present 
section. We begin with the following proposition, which serves as the first step 
in this direction.

\bigskip

\begin{proposition}
    Let $(\m,g)$ be an $N$-dimensional Cartan–Hadamard model manifold with $N \geq 3$. Suppose $u$ is a measurable function on $\m$. Then  
    \begin{align}\label{comparison}
        \|u\|_{L^{2^\star,2}(\m)}\leq \left(\frac{\uv}{N}\right)^{-\frac{1}{N}}\left(\int_{\m}\frac{(u^\sharp)^2}{r^2}\dvg\right)^{\frac{1}{2}},
    \end{align}
    where $u^\sharp$ and $L^{2^\star,2}(\m)$ are defined in \eqref{sharp} and \eqref{org-quasinorm}, respectively.
\end{proposition}
\begin{proof}
       Since $u^\sharp$ is radially decreasing, we write, by abusing notation, $u^\sharp(x)=u^\sharp(r):=f(r)$, where $r=\mathrm{dist}_g(x,x_0)$. Now, by definition, $u^\star(V(r))=f(r)$ for $r>0$ and recall $V(r)$ is defined in \eqref{vol}. Equivalently, we write $(u^\sharp)^\star(s)=f(G(s))$ for $s>0$, where $s=V(r)$ and $r=G(s)$. Now, doing a change of variable $s=V(r)$ and using $V'(r)=\omega_N\,(\psi(r))^{N-1}$, we have
    \begin{align*}
        \|u^\sharp\|^2_{L^{2^\star,2}(\m)}=\int_0^{\infty}s^{\frac{2}{2^\star}-1}f(G(s))^2\ds=\int_0^{\infty}V(r)^{\frac{2}{2^\star}-1}(f(r))^2V'(r)\dr=\int_{\m}V(r)^{-\frac{2}{N}}(u^\sharp)^2\dvg.
    \end{align*}
Now using \eqref{sturm}, we deduce $V(r)\geq \frac{\uv}{N}\,r^{N}$, and using this in the above, along with the rearrangement invariance property, i.e., $ \|u\|^2_{L^{2^\star,2}(\m)}=\|u^\sharp\|^2_{L^{2^\star,2}(\m)}$,  we establish \eqref{comparison}.
\end{proof}

\medskip

 Before presenting the next lemma, we recall a key result of Kufner and Persson 
\cite[Theorem~6.14, (iv)]{kuf-per-03}, which will be essential for our argument. 
In particular, for our purpose, it can be stated as follows:
\begin{lemma}[\cite{kuf-per-03}]\label{sup-condition}
Let $N\geq 3$ and suppose $\vartheta_1$ and $\vartheta_2$ be two non-negative Borel measures on $(0,\infty)$. Then
\begin{equation}\label{new inequality}
        \left(\int_0^{\infty}f(t)^{\frac{N+2}{N-2}}\:\mathrm{d}\vartheta_1(t)\right)^{\frac{N-2}{N+2}}\leq C\int_0^{\infty}f(t)\:\mathrm{d}\vartheta_2(t)
\end{equation}
holds for all non-negative, non-increasing Borel measurable functions $f$ on $(0,\infty)$ if and only if
\begin{equation*}
        B:= \sup_{r>0}\frac{\left(\int_0^r\:\mathrm{d}\vartheta_1(t)\right)^{\frac{N-2}{N+2}}}{\left(\int_0^r\:\mathrm{d}\vartheta_2(t)\right)}<+\infty,
\end{equation*}
where $C$ is a finite positive constant independent of $f$ and the r.h.s. of \eqref{new inequality} is assumed to be finite. Furthermore, $C=B$ is the sharp constant.
\end{lemma}

\medskip

Recall that, by a standard result from rearrangement theory, the Hardy–Littlewood inequality (see \cite[Theorem~2.2]{BS}) gives
\begin{equation}\label{Hardy Littlewood}
    \int_{\m}\frac{u^2}{r^2}\,\mathrm{d}v_g \leq \int_{\m}\frac{(u^\sharp)^2}{r^2}\,\mathrm{d}v_g,
\end{equation}
for any admissible function $u$ defined on the Cartan–Hadamard model manifold. We now perform a stability analysis of \eqref{Hardy Littlewood} by cleverly applying \cite[Theorem~1.3 and Remark~1.4]{cf-jlms}. It is worth mentioning that the corresponding result in the flat manifold was first established by Cianchi-Ferone in \cite[Lemma~2.1]{cf-aihp}. Here, we extend that analysis to general Cartan–Hadamard model manifolds. Here is our result below.

\medskip

\begin{lemma}\label{main-lem}
   Let $(\m,g)$ be an $N$-dimensional Cartan–Hadamard model manifold with $N \geq 3$. Then a positive constant $C = C(N)$ exists such that
    \begin{equation}\label{HL stability}
        \int_{\m} \frac{u^2}{r^2}\dvg+C\left(\int_{\m}\frac{(u^\sharp)^2}{r^2}\dvg\right)^{-\frac{N+2}{N-2}}\left(\int_{\m}|u(x)-u^\sharp(x)|^{2^\star}\dvg\right)^2\leq \int_{\m}\frac{(u^\sharp)^2}{r^2}\dvg,
    \end{equation}
    holds for every non-negative measurable function $u\in D^{1,2}(\m)$ and making the r.h.s. of \eqref{HL stability} finite.
\end{lemma}
\begin{proof}
Suppose $g:\m \to [0,\infty)$ is a radially strictly decreasing function that decays to zero at infinity, and assume that the function $\theta:(0,\infty)\to [0,\infty)$ defined by
\begin{equation*}
    \theta(s) = {\text{ess}\sup}_{r\in (0,s)} \left(-\frac{1}{{g^\star}'(r)}\right), \qquad \text{ for all }\quad s>0,
\end{equation*}
is finite, locally absolutely continuous, and satisfies $\lim_{s\to 0^+} \theta(s) = 0$.

Let $f:\m \to [0,\infty)$ be any function decaying to zero at infinity, such that the quasinorm
\begin{equation*}
    \|f\|_{\Lambda^q} = \left(\int_0^{\infty} f^\star(s)^q \, \theta'(s) \ds\right)^{\frac{1}{q}}
\end{equation*}
is finite for some $q \in [1,\infty)$. Then, thanks to the Hardy–Littlewood inequality with a remainder term \cite[Theorem~1.3 and Remark~1.4]{cf-jlms} and after applying it for a particular case, there exists a uniform constant $C>0$ such that
\begin{equation}\label{HL first stability}
    \int_{\m} f(x)g(x)\,\mathrm{d}v_g 
    + \frac{1}{Cq}\,\|f\|_{\Lambda^q}^{-q}\,\,\|f - f^\sharp\|_{L^{\frac{q+1}{2}}(\m)}^{q+1}
    \leq \int_{\m} f^\sharp(x)g(x)\,\mathrm{d}v_g.
\end{equation}

Now, we will choose $f(x)= u(x)^2$, radially strictly decreasing function $g(x)=r^{-2}$ and an integrability exponent $q=\frac{N+2}{N-2}$. Then we need to verify that all the hypotheses hold in our setting.

First, by direct calculation, we have $g^\star(s)=\left(G(s)\right)^{-2}$ for $s>0$. First, using the chain rule, we calculate
\begin{align*}
   -\frac{1}{{g^\star}'(s)}=\frac{1}{2}\frac{\left(G(s)\right)^3}{(G)'(s)}=\frac{\uv}{2}(G(s))^3(\psi(G(s)))^{N-1}. 
\end{align*}
Then, using \eqref{ins-assump}, we observe $(G(s))^3(\psi(G(s)))^{N-1}$ is an increasing function and so 
\begin{align*}
    \theta(s)=\frac{\uv}{2}(G(s))^3(\psi(G(s)))^{N-1} \qquad \text{ for all }\quad s>0.
\end{align*}
Moreover, it is finite and locally absolutely continuous. Also, using $\psi(0)=0$ and $G(0)=0$ we deduce $\lim_{s\to 0^+}\theta(s)=0$.
Next, we shall verify that
\begin{align}\label{work}
    \|f\|_{\Lambda^q} \leq C\int_{\m}\frac{{u^\sharp}^2}{r^2}\dvg,
\end{align}
where $C$ is some positive generic constant. Now, using change of variable $s=V(r)$, we deduce
\begin{align*}
    \|f\|_{\Lambda^q}^q=\int_0^{\infty} f^\star(s)^q \, \theta'(s) \ds=\int_0^{\infty}\left((u^\sharp(r))^2\right)^{q}\left(\frac{\uv}{2}r^3\psi(r)^{N-1}\right)'\dr,
\end{align*}
and
\begin{align*}
   \int_{\m}\frac{{u^\sharp}^2}{r^2}\dvg=  \int_0^{\infty}(u^\sharp(r))^2\left(\uv \frac{\psi(r)^{N-1}}{r^2}\right)\dr=\int_0^{\infty}(u^\sharp(r))^2\left(\int_0^r \uv \frac{\psi(t)^{N-1}}{t^2}\dt\right)'\dr.
\end{align*}
Now, thanks to Lemma~\ref{sup-condition} and Lemma~\ref{complete-lem}, \eqref{work} follows.

Therefore, for these specific choices of $f$, $g$ , $q$ and applying it in \eqref{HL first stability}, along with \eqref{work} and $(u^\sharp)^2=(u^2)^\sharp$, we obtain
 \begin{equation}\label{finalchange}
        \int_{\m}\frac{u(x)^2}{r^2}\dvg+C\left(\int_{\m}\frac{(u^\sharp)^2}{r^2}\dvg\right)^{-\frac{N+2}{N-2}}\|u^2-(u^\sharp)^2\|_{L^{\frac{N}{N-2}}(\m)}^{2^\star}\leq\int_{\m}\frac{(u^\sharp)^2}{r^2}\dvg.
\end{equation}
   
Since $|s-r|^2\leq |s^2-r^2|$ for every $r,s\geq0$, we have
\begin{equation}\label{first change}
            \|u-u^\sharp\|^{\frac{4N}{N-2}}_{L^{2^\star}(\m)}\leq \|u^2-(u^\sharp)^2\|_{L^{\frac{N}{N-2}}(\m)}^{2^\star}
\end{equation}
Combining \eqref{first change} in \eqref{finalchange}, the desired estimate \eqref{HL stability} follows.
\end{proof}

Having developed the required preliminary results, we are now equipped to formulate 
and prove the main theorem of this section, which presents the core conclusion 
of our stability analysis.

\begin{theorem}\label{main-th-hardy-cf-stab}
Let $(\m,g)$ be an $N$-dimensional Cartan--Hadamard model manifold with $N \geq 3$. 
Consider any function $u \in D^{1,2}(\m)$ and recall the distance from $u$ to the 
family of extremal profiles $\{\mathcal{U}_a\}_{a\in \mathbb{R}}$, defined by
\begin{equation*}
    d_{\m}(u) := \inf_{a\in\mathbb{R}}
    \frac{\|u\;-\;\mathcal{U}_a\|_{\widetilde{L}^{2^\star,\infty}(\m)}}{\|u\|_{D^{1,2}(\m)}}.
\end{equation*}
 Then there exists a constant $C = C(N) > 0$ such that the following \emph{quantitative} 
Hardy inequality holds:
\begin{equation*}
    \left(\frac{N-2}{2}\right)^2 
    \int_{\m}\frac{u(x)^2}{r^2}\, \dvg
    \left( 1 + C\,(d_{\m}(u))^{\frac{4N}{N-2}} \right)
    \leq 
    \int_{\m} |\nabla_g u|^2\, \dvg.
\end{equation*}
\end{theorem}

\medskip

\begin{proof}
We begin the proof by assuming 
\begin{equation}\label{first assumption}
        u\geq 0,
    \end{equation}
    and
        \begin{equation}\label{second assumption}
         \int_{\m}{\frac{u^\sharp(x)^2}{r^2}\dvg}=1.
     \end{equation}
For notational simplicity, we assume $C_H:=\frac{(N-2)^2}{4}$. Let us recall the definition of $\eu(u)$ from \eqref{energy-hardy} , and we rewrite it as follows
\begin{align}\label{3f}
    \eu (u) =\left(\int_{\m}|\gradg u|^2\dvg-\int_{\m}|\gradg u^\sharp|^2\dvg\right)  
          +\left(\int_{\m}|\gradg u^\sharp|^2\dvg-C_H\right)+C_H\left(1-\int_{\m}\frac{u(x)^2}{r^2}\dvg\right).
\end{align}

The first factor of the above integral is nonnegative by the P\'olya–Szeg\"o inequality. Again, using the Hardy-Littlewood inequality and \eqref{second assumption}, we have that the third factor in \eqref{3f} is non-negative. For the second factor, from the assumed condition $\int_{\m}\frac{u^\sharp(x)^2}{r^2}\dvg=1$, applying \eqref{hardy} for $u^\sharp$, P\'olya–Szeg\"o inequality, and Hardy-Littlewood inequality, we deduce
\begin{align*}
    \eu(u)\geq \int_{\m}|\gradg u^\sharp|^2\dvg-C_H= \eu(u^\sharp)\geq 0.
\end{align*}

Again, by exploiting \eqref{second assumption} together with Proposition~\ref{1D hardy-prop} for the function $u^\sharp$, we deduce
\begin{align*}
 \eu(u^\sharp)= \uv^2\int_0^{\infty} |{u^\star}'(s)|^2\left(\psi \left(G(s)\right)\right)^{2(N-1)}\ds- C_H\int_0^{\infty}\frac{u^\star(s)^2}{G(s)^2}\ds. 
\end{align*}
Therefore, applying the above, \eqref{second assumption} in Proposition~\ref{1stabi} for the function $u^\sharp$, and then simplifying \eqref{L to Ltilde}, we receive
\begin{equation}\label{fnl-step1}
\inf_{a\geq 0}\left\|u^\star(s)-a\left(G(s)\right)^{\frac{2-N}{2}}\right\|_{\tilde{L}_{\m}^{2^\star,\infty}(0,\infty)}\leq C \, \eu(u)^{\frac{1}{4}},
\end{equation}
where $C$ is some generic positive absolute constant.

 Moreover, applying our main Lemma~\ref{main-lem}, \eqref{3f}, and \eqref{second assumption}, we get
 \begin{equation*}
          \left(\int_{\m}\left|u(x)-u^\sharp(x)\right|^{2^\star}\dvg \right)^2\leq C\,\eu (u).
\end{equation*}
Then, applying the chain of embeddings \eqref{fnl-emb} in the above, we have
 \begin{equation}\label{fnl-step2}
\|u-u^\sharp\|_{\tilde{L}^{2^\star ,\infty}(\m)}\leq C\, \eu (u)^{\frac{N-2}{4N}}.
\end{equation}

Now, recall the extremizer $\{\mathcal{U}_a\}$ in \eqref{extremizer} and combining \eqref{fnl-step1} and \eqref{fnl-step2}, we obtain, observing that $u^\sharp \;- \; \mathcal{U}_a$ and $u^\star(s)-a(G(s))^{\frac{2-N}{2}}$ are equimeasurable,
\begin{align}\label{yeld}
      \inf_{a\geq 0} \|u-\mathcal{U}_a\|_{\tilde{L}^{2^\star,\infty}(\m)}
          & \leq C\left(\|u-u^\sharp\|_{\tilde{L}^{2^\star,\infty}(\m)}+\inf_{a\geq 0}\|u^\sharp-\mathcal{U}_a\|_{\tilde{L}^{2^\star,\infty}(\m)}\right)\nonumber \\ &
          =C\left(\|u-u^\sharp\|_{\tilde{L}^{2^\star,\infty}(\m)}+\inf_{a\geq 0}\left\|u^\star(s)-a\left(G(s)\right)^{\frac{2-N}{2}}\right\|_{\tilde{L}_{\m}^{2^\star,\infty}(0,\infty)}\right)\nonumber\ \\ & \leq 
           C\left(\eu(u)^{\frac{N-2}{4N}}+\eu(u)^{\frac{1}{4}}\right).
\end{align}

If $\eu(u)\leq 1$, then using $\frac{N-2}{4N}\leq \frac{1}{4}$ in \eqref{yeld} we get 
 \begin{equation}\label{yeld-1}
          \inf_{a\geq 0} \|u-\mathcal{U}_a\|_{\tilde{L}^{2^\star,\infty}(\m)}\leq C\eu(u)^{\frac{N-2}{4N}}.
\end{equation}
Whenever $\eu(u)> 1$, then using \eqref{emb}, \cite[Prop.~4.2, p.~217]{BS}, \eqref{comparison}, we have
\begin{multline}\label{yeld-2}
     C\inf_{a\geq 0}\|u-\mathcal{U}_a\|_{\tilde{L}^{2^\star,\infty}(\m)} \leq C \|u\|_{\tilde{L}^{2^\star,\infty}(\m)} \leq C \|u\|_{L^{2^\star,\infty}(\m)}\\ \leq  C \|u\|_{L^{2^\star,2}(\m)}\leq \int_{\m}\frac{(u^\sharp)^2}{r^2}\dvg=1<\eu(u).
\end{multline}
Therefore, combining \eqref{yeld-1} and \eqref{yeld-2}, we obtain
 \begin{equation}\label{derive-1}
          \left(\frac{N-2}{2}\right)^2\int_{\m}\frac{u(x)^2}{r^2}\dvg+C\inf_{a\geq 0}\|u-\mathcal{U}_a\|^{\frac{4N}{N-2}}_{\tilde{L}^{2^\star,\infty}(\m)}\leq\int_{\m}|\gradg u|^2\dvg,
\end{equation}
for some positive constant $C=C(N)$ under the two assumptions \eqref{first assumption}, and \eqref{second assumption} specified at the very beginning. First, we relax the assumption \eqref{second assumption}, by considering $u\left(\int_{\m}\frac{(u^\sharp)^2}{r^2}\dvg\right)^{-\frac{1}{2}}$ in the place of $u$ in \eqref{derive-1} and we derive
\begin{align}\label{derive-2}
   \inf_{a\geq 0}\|u-\mathcal{U}_a\|^{\frac{4N}{N-2}}_{\tilde{L}^{2^\star,\infty}(\m)}\leq C \left(\int_{\m}\frac{(u^\sharp)^2}{r^2}\dvg\right)^{2^\star-1}\left(\int_{\m}|\gradg u|^2\dvg- \left(\frac{N-2}{2}\right)^2\int_{\m}\frac{u^2}{r^2}\dvg\right).
\end{align}

\medskip 

Finally, we relax the assumption \eqref{first assumption} by considering $u_+=\frac{|u|+u}{2}$ and $u_-=\frac{|u|-u}{2}$ in \eqref{derive-2}. Here we follow exactly the same approach as in the proof  of \cite[Theorem~1.1, Eqn.~2.50, 2.51]{cf-aihp}, and we deduce
\begin{align*}
    \inf_{a\in \mathbb{R}}\|u-&v_a\|_{\tilde{L}^{2^\star,\infty}(\m)} \leq  C\left(\int_{\m}\frac{(u^\sharp)^2}{r^2}\dvg\right)^{\frac{1}{2}}\left(\int_{\m}\frac{u^2}{r^2}\dvg\right)^{-\frac{N-2}{4N}}\\& \qquad\qquad \times \sum_{\pm}\left(\int_{\m}|\gradg u_{\pm}|^2\dvg- \left(\frac{N-2}{2}\right)^2\int_{\m}\frac{u_{\pm}^2}{r^2}\dvg\right)^{\frac{N-2}{4N}}\\& \leq C\|u\|_{D^{1,2}(\m)}\left(\int_{\m}\frac{u^2}{r^2}\dvg\right)^{-\frac{N-2}{4N}}\left(\int_{\m}|\gradg u|^2\dvg- \left(\frac{N-2}{2}\right)^2\int_{\m}\frac{u^2}{r^2}\dvg\right)^{\frac{N-2}{4N}},
\end{align*}
where in the end we used $\left(\int_{\m}\frac{(u^\sharp)^2}{r^2}\dvg\right)^{\frac{1}{2}}\leq C \|u\|_{D^{1,2}(\m)}$ which follows by application of the Hardy inequality and the Pólya-Szegő inequality. This completes the thesis.
\end{proof}

\begin{remark}
{\rm 
In particular, Theorem~\ref{main-th-hardy-cf-stab} shows that when 
$d_{\m}(u)$ is sufficiently small---meaning that $u$ is close to the family of 
extremal functions $\{\mathcal{U}_a\}$---the resulting inequality yields a genuine 
improvement over the classical sharp Hardy inequality. More precisely, the deficit in the inequality controls the distance of $u$ from the extremal family $\{\mathcal{U}_a\}$ in a quantitative way.
}
\end{remark}

\begin{remark}
{\rm
 We remark that Theorem~\ref{main-th-hardy-cf-stab} exhibits a slight variation from 
\cite[Theorem~1.1]{cf-aihp}. In the latter, the distance from the extremal family 
is measured using the norm $\|u\|_{L^{2^\star,2}(\mathbb{R}^N)}$, whereas in our 
framework we employ the Dirichlet norm $\|u\|_{D^{1,2}(\m)}$ in the definition of 
the distance function \eqref{hardy stability distance}. This modification is 
dictated by the underlying geometry of the manifold: in our case, the quantity 
\[
    \left(\int_{\m}\frac{(u^\sharp)^2}{r^2}\dvg\right)^{1/2}
\]
is naturally controlled by the Dirichlet norm, rather than by a 
$\widetilde{L}^{2^\star,2}(\m)$-type norm. Indeed, if an estimate of the form
\[
    \left(\int_{\m}\frac{(u^\sharp)^2}{r^2}\dvg\right)^{1/2} 
    \leq C\,\|u\|_{\widetilde{L}^{2^\star,2}(\m)}
\]
were valid on $\m$, then the geometry would necessarily coincide with the Euclidean 
one (see Remark~\ref{only-euclid}). Hence, our choice of norm in the stability 
distance reflects an essential geometric feature of the non-Euclidean model 
manifold.

    }
\end{remark}
\medskip 

\section{Optimal Sobolev-Lorentz embedding in a general Cartan-Hadamard model manifold}\label{sec-SL}

It is well known that the classical Sobolev inequality admits refinements when formulated within the framework of Lorentz spaces in the Euclidean setting. These strengthened versions yield sharper embeddings and better capture the fine behaviour of functions near the critical Sobolev exponent. Significant progress in this direction has been made in the foundational works of Peetre \cite{pet} and Alvino \cite{alv}. In particular, the refinements mentioned  lead to the following embedding results:

\begin{align}\label{pert-alvi}
  D^{1,2}(\mathbb{R}^N)
  \hookrightarrow 
  L^{2^\star, 2}(\mathbb{R}^N)
  \hookrightarrow 
  L^{2^\star}(\mathbb{R}^N).
\end{align}

The significance of the embedding~\eqref{pert-alvi} stems from the fact that it is 
\emph{optimal} within the class of rearrangement invariant spaces (see~\cite{ekp}). 
In particular, among all possible rearrangement invariant quasinorms, the Lorentz 
quasinorm $\|\cdot\|_{L^{2^\star,2}(\mathbb{R}^N)}$ is the finest one that can appear 
on the left-hand side of the inequality. Equivalently, no strictly stronger 
rearrangement invariant space may replace $L^{2^\star,2}(\mathbb{R}^N)$ while still 
preserving the validity of the following sharp embedding:
\begin{align}\label{pick-best}
    \|u\|_{L^{2^\star,2}(\mathbb{R}^N)}
    \leq S_{N,2^\star,2}\, \|u\|_{D^{1,2}(\mathbb{R}^N)}.
\end{align}
The optimal constant in~\eqref{pick-best} is explicitly given by
\begin{align}\label{alv-best}
    S_{N,2^\star,2}
    = \frac{2}{N-2}\,
      \frac{\big[\Gamma\left(1+\frac{N}{2}\right)\big]^{\frac{1}{N}}}{\sqrt{\pi}},
\end{align}
was computed by Alvino in his seminal work \cite{alv}; see also the related contributions in \cite{crt, crt2}. This result highlights the sharpness of the Lorentz refinement of the classical Sobolev inequality.

\medskip 

Using \eqref{comparison} together with the Hardy inequality \eqref{hardy}, the P\'olya--Szeg\H{o} principle \eqref{ps}, and following the argument of \cite[Prop.~4.2, p.~217]{BS}, we can establish \eqref{pert-alvi} on Cartan--Hadamard model manifolds. This Sobolev--Lorentz embedding can be written as
\begin{align}\label{slmm}
    D^{1,2}(\m) \hookrightarrow L^{2^\star,2}(\m) \hookrightarrow L^{2^\star}(\m).
\end{align}

\medskip

To the best of our knowledge, the analogous result was previously known only in the 
hyperbolic setting, where it was obtained by Nguyen in~\cite[Theorem~1.2]{VHN} using 
a different technique. For general Cartan--Hadamard model manifolds, this result appears 
to be new. In this section, we show that there exists no proper rearrangement-invariant 
space lying strictly between $D^{1,2}(\m)$ and $L^{2^\star}(\m)$; in this sense, the embedding 
\eqref{slmm} is optimal. Moreover, we compute
\begin{align*}
    \sup_{\,u \in D^{1,2}(\m),\, \|u\|_{D^{1,2}(\m)} \neq 0}
    \frac{\|u\|_{L^{2^\star,2}(\m)}}{\|u\|_{D^{1,2}(\m)}},
\end{align*}
which yields the sharp constant for the embedding. Finally, we prove that this optimal 
constant is never achieved.

The transformation \eqref{transformation}, introduced earlier in Section~1, serves as the central ingredient in extending the result to the present geometric setting. Its structural properties allow us to transfer the problem to a setting where rearrangement techniques and sharp estimates become applicable. Before formulating the main theorem, we first establish several auxiliary lemmas that will play a key role in the proof.

Let $(\m,g)$ be an $N$-dimensional Cartan–Hadamard model manifold with $N \geq 3$. We define a map $\st:\m \to \rn$ by
\begin{align*}
    \st(x)=\st(r,\theta)=(\varrho, \theta) \qquad \text{ for all } x\in \m,
\end{align*}
where $x=(r,\theta)$ and $r=\text{dist}(x_0,x)$, $\theta\in\sn$ and $\varrho$ satisfy the \eqref{transformation}, i.e.,
\begin{align}\label{trans}
\uv \int_0^{r}(\psi(t))^{N-1}\dt = \uv \int_0^{\varrho}t^{N-1}\dt.
\end{align}

Now by differentiation \eqref{trans}, we observe that for $r>0$, there holds
\begin{align}\label{dvol-pre}
  \psi(r)^{N-1}\dr=\varrho^{N-1}\drho.
\end{align}

Let $\mathcal{M}(\cdot)$ denote the set of all real-valued measurable functions on the underlying space. We now define $\scal: \mathcal{M}(\m) \to \mathcal{M}(\rn)$ by
\begin{align}\label{transport}
\scal(f):= f\circ\st^{-1} \qquad \text{ for all } f\in \mathcal{M}(\m),
\end{align}
where $\st^{-1}$ exists as $\st$ is a homeomorphism between $\m$ and $\rn$ via polar coordinate and our transformation defined in \eqref{trans}. Moreover, it is easy to see that $\scal$ is a bijection and it preserves the linear structure.

Now, we are ready to state our required lemmas.
\begin{lemma}
 Let $f\in \mathcal{M}(\m)$. Then $\scal(f)$ defined in \eqref{transport} and $f$ are equidistributed.   Moreover, $\scal\big|_{L^{2^\star,2}(\m)}:L^{2^\star,2}(\m)\to L^{2^\star,2}(\rn)$ is an isometry, i.e,
 \begin{align}\label{iso-lem}
     \|f\|_{L^{2^\star,2}(\m)}=\|\scal(f)\|_{L^{2^\star,2}(\rn)} \qquad \text{ for all } f\in L^{2^\star,2}(\m).
 \end{align}
\end{lemma}
\begin{proof}
  We start with $\lambda>0$ and we observe 
  \begin{align*}
      \st\left(\{x \in \m: |f(x)|> \lambda\}\right)=\{y \in \rn: |(\scal(f))(y)|> \lambda\}.
  \end{align*}
  Therefore, using the volume preserving property of $\st$, we deduce
   \begin{align}\label{equi-dis}
     |\{x \in \m: |f(x)|> \lambda\}|=\mathcal{L}^N(\{y \in \rn: |(\scal(f))(y)|> \lambda\}).
  \end{align}
The above implies $\scal(f)$ and $f$ are equidistributed, and by definition \eqref{org-quasinorm}, the quasinorm is preserved.
\end{proof}
\begin{lemma}
  Suppose $f\in D_{\text{rad}}^{1,2}(\m)$  then $\scal(f)\in D_{\text{rad}}^{1,2}(\rn)$. Moreover, it holds 
  \begin{align}\label{eq-sobo-rad-lem}
      \|\scal(f)\|_{D_{\text{rad}}^{1,2}(\rn)}\leq  \|f\|_{D_{\text{rad}}^{1,2}(\m)} \qquad \text{ for all }\quad f\in D_{\text{rad}}^{1,2}(\m).
  \end{align}
\end{lemma}
\begin{proof}
    For notational simplicity, we denote $v=\scal(f)$. It is immediate that as $f$ is radial, then $v$ is also radial and by definition $f(r)=v(\varrho(r))$, where $r$ and $\varrho$ are connected by \eqref{trans}. Now, writing from $\varrho$ to $r$ and using \eqref{dvol-pre}, we notice that 
    \begin{align}\label{sobo-rad-lem-1}
        \|v\|_{D_{\text{rad}}^{1,2}(\rn)}^2=\uv\int_0^\infty\left( \frac{\partial v }{\partial \varrho}(\varrho)\right)^2\varrho^{N-1}\drho=\uv\int_0^\infty\left( \frac{\partial f }{\partial r}(r)\right)^2(h(r))^{2(N-1)}(\psi(r))^{N-1}\dr,
    \end{align}
    where
    \begin{align}\label{hcall}
      h(r)=\frac{\varrho(r)}{\psi(r)}  \qquad \text{ for all } \quad r>0.
    \end{align}
    Now, by \eqref{trans} and using \eqref{sturm}, we obtain
    \begin{align*}
        \varrho(r)^N=N\int_0^r\psi(t)^{N-1}\dt\leq N\int_0^r\left(\frac{t}{r}\psi(r)\right)^{N-1}\dt=r\psi(r)^{N-1}\leq \psi(r)^N.
    \end{align*}
    Then using this we deduce $h(r)\leq 1$ for every $r>0$ and hence using this in \eqref{sobo-rad-lem-1} and writing back to the expression of $ \|f\|_{D_{\text{rad}}^{1,2}(\m)}^2$ via polar coordinate we complete the proof.
\end{proof}

\medskip 

\begin{remark}
{\rm 
 Let $(\m,g)$ be an $N$-dimensional Cartan–Hadamard model manifold with $N \geq 3$. From \eqref{trans} and using the first part of \eqref{sturm}, we notice that to create a fixed amount of volume, we require more Euclidean radius $\varrho(r)$ than the underlying Cartan-Hadamard model manifold radius $r$. Moreover, it also follows from the above 
 
 $$
 \dfrac{r}{\psi(r)} \leq 1  \implies  \dfrac{\varrho(r)}{\psi(r)} \leq 1.
 $$
}
where $r$ and $\varrho(r)$ are related through the \eqref{trans}.
\end{remark}

\begin{lemma}\label{img-ri}
Suppose $X(\m)$ is a rearrangement invariant function space. Then $\scal(X(\m))$ is a rearrangement invariant function space on $\rn$.
\end{lemma}
\begin{proof}
 Let us denote $Y(\rn):=\scal(X(\m))=\{\scal(f): f \in X(\m)\}$. Then we define
 \begin{align}\label{ri-dfn}
    \|u\|_{Y(\rn)}:=\|\scal^{-1}(u)\|_{X(\m)} \qquad \text{ for all } \quad  u\in Y(\rn).
\end{align}
Thanks to the linear structure of $\scal^{-1}$, which implies $Y(\rn)$ is a complete normed linear space. Now suppose $u$ and $v$ are two functions on $Y(\rn)$ with $u^\star=v^\star$. Also by definition there exists $f$ and $g$ in $X(\m)$ such that $u=\scal(f)$ and $v=\scal(g)$. Then, using \eqref{equi-dis}, we can see that $f^\star=g^\star$ on $\m$. As given $X(\m)$ is a rearrangement invariant function space, hence
\begin{align*}
    \|f\|_{X(\m)}=\|g\|_{X(\m)}.
\end{align*}
Finally, in above using $f=\scal^{-1}(u)$ and $g=\scal^{-1}(v)$ and the definition \eqref{ri-dfn}, we deduce
\begin{align*}
   \|u\|_{Y(\rn)}=\|v\|_{Y(\rn)}.
\end{align*}
\end{proof}

\begin{lemma}
    Suppose $X(\m)$ is a norm linear space and  $X(\m)\hookrightarrow L^{2^\star,2}(\m)$ is a proper embedding. Then  $\scal(X(\m))\hookrightarrow L^{2^\star,2}(\rn)$ is also a proper embedding.
\end{lemma}
\begin{proof}
    We notice that by assigning definition \eqref{ri-dfn}, the space $\scal(X(\m))$ is a complete normed linear space. Suppose $u=\scal(f)\in \scal(X(\m))$ for some $f\in X(\m)$. Then, using \eqref{iso-lem}, embedding, and \eqref{ri-dfn}, we have
    \begin{align*}
        \|u\|_{L^{2^\star,2}(\rn)}=\|\scal(f)\|_{L^{2^\star,2}(\rn)}=\|f\|_{L^{2^\star,2}(\m)}\leq C \|f\|_{X(\m)}=C\|u\|_{\scal(X(\m))}.
    \end{align*}
    This proves that $\scal(X(\m))\hookrightarrow L^{2^\star,2}(\rn)$ and it is proper due to the properness of $X(\m)\hookrightarrow L^{2^\star,2}(\m)$.
\end{proof}

\begin{lemma}\label{rad-lem-nh}
    The embedding $D^{1,2}(\m)\hookrightarrow L^{2^\star,2}(\m)$ is sharp in the sense that there is no rearrangement invariant function space $X(\m)$ for which
    \begin{align}\label{ifpossible}
        D^{1,2}(\m)\hookrightarrow X(\m)\hookrightarrow L^{2^\star,2}(\m)
    \end{align}
    holds and $ D^{1,2}(\m)\hookrightarrow X(\m)$ and $X(\m)\hookrightarrow L^{2^\star,2}(\m)$ are proper embedding.
\end{lemma}
\begin{proof}
    If possible, some rearrangement invariant function space $X(\m)$ exists and \eqref{ifpossible} holds. It immediately gives
      \begin{align}\label{rifpossible}
        D_{\text{rad}}^{1,2}(\m)\hookrightarrow X(\m)\hookrightarrow L^{2^\star,2}(\m).
    \end{align}  

    Now using first part of \eqref{rifpossible} and density argument, there exists $u\in C^{\infty}_{c,\text{rad}} (\rn)$ such that its $D_{\text{rad}}^{1,2}(\rn)$ norm is finite whereas the $\scal(X(\m))$ norm is infinite. Indeed, if $D^{1,2}_{\text{rad}}(\rn)$ is embedded inside $\scal(X(M))$ then via Pólya-Szegő inequality and rearrangement invariance property of $\scal(X(\m))$ (refer Lemma~\ref{img-ri})  we will have $D^{1,2}(\rn) \hookrightarrow \scal(X(\m))$ which contradicts the optimality of the embedding \eqref{pert-alvi}.

   Then writing $f=\scal^{-1}(u)$ we have $\|f\|_{X(\m)}=\|u\|_{\scal(X(\m))}=+\infty$. Moreover, a similar computation as in  \eqref{sobo-rad-lem-1}, we deduce
\begin{align*}
        \|f\|_{D^{1,2}(\m)}^2 =\int_{\rn} |\nabla u(\varrho)|^2 (h(r(\varrho)))^{-2(N-1)}\dx<+\infty,
\end{align*}
where we used the facts $u$ is compactly supported and $h$ defined earlier in \eqref{hcall}. Also, we used the fact that $h$ is a continuous function and 
\begin{align}\label{h-lomit}
    \lim_{\varrho\to 0^+} h(r(\varrho))=1
\end{align}
which follows by exploiting first part of \eqref{sturm} and \eqref{psi}.  Thus, we have produced one radial function $f$ which is in $D^{1,2}(\m)$ but not in $X(\m)$, and this contradicts \eqref{rifpossible}. Hence the embedding $D^{1,2}(\m)\hookrightarrow L^{2^\star,2}(\m)$ is sharp. 
\end{proof}
\medskip

We are now prepared to present the optimal Sobolev–Lorentz embedding theorem on Cartan–Hadamard model manifolds. The groundwork established in the preceding lemmas allows us to rigorously formulate the result and proceed with its proof.

\medskip 

\begin{theorem}
Let $(\m,g)$ be a Cartan-Hadamard model manifold with dimension $N\geq 3$. Then $D^{1,2}(\m)$ is continuously embedded in $L^{2^\star,2}(\m)$ which is the smallest rearrangement-invariant space containing $D^{1,2}(\m)$. Furthermore, the sharp constant of this embedding is explicitly given by 
\begin{align*}
    \|u\|_{L^{2^\star,2}(\m)} \leq S_{N,2^\star,2}\|u\|_{D^{1,2}(\m)}
\end{align*}
where $S_{N,2^\star,2}$ is defined in \eqref{alv-best} and this constant is never attained in $D^{1,2}(\m)$.
\end{theorem}

\begin{proof}
The optimality of the embedding \eqref{slmm} can be seen from  Lemma~\ref{rad-lem-nh}. Now we will show that the above constant $S_{N,2^\star,2}$ is optimal and the constant coincides with the best constant as in the Euclidean space. We will prove
\begin{align}\label{f-1}
    \sup_{u\in D^{1,2}(\m),\, \|u\|_{ D^{1,2}(\m)}\neq 0}
    \frac{\|u\|_{ L^{2^\star,2}(\m)}}{\|u\|_{ D^{1,2}(\m)}}=S_{N,2^\star,2}.
\end{align}
By the Pólya–Szeg\"o inequality, it is enough to consider the Rayleigh quotient in \eqref{f-1} over all radial functions. Let $f\in D^{1,2}_{\text{rad}}(\m)$ and denote $u=S(f)$. Now using \eqref{eq-sobo-rad-lem}, \eqref{ri-dfn}, \eqref{pick-best} and taking supremum, we deduce
\begin{align*}
    \sup_{u\in D^{1,2}(\m),\, \|u\|_{ D^{1,2}(\m)}\neq 0}
    \frac{\|u\|_{ L^{2^\star,2}(\m)}}{\|u\|_{ D^{1,2}(\m)}}\leq S_{N,2^\star,2}.
\end{align*}

For the reverse inequality, consider $u\in C^{\infty}_{c,\,\text{rad}}(\rn)$. For every $k\in\mathbb{N}$, we define $f_k$ on $\m$ by
\begin{align*}
    f_k:=\scal^{-1}(u_k) \qquad \text{ where } \quad u_k(y)= u(ky) \quad \text{ for every } y\in\rn.
\end{align*}

Therefore, we will be done if we establish
\begin{align}\label{fnl-claim-petri}
    \lim_{k \to \infty}\frac{\|f_k\|_{L^{2^\star,2}(\m)}}{\|f_k\|_{D^{1,2}(\m)}}=\frac{\|u\|_{_{L^{2^\star,2}(\rn)}}}{\|u\|_{D^{1,2}(\rn)}}.
\end{align}

Next, by definition $f_k(r)=\scal^{-1}u_k(r)=u_k(\varrho(r))=u(k\varrho(r))$ and denote $H(\varrho)=(h(r(\varrho)))^{-2(N-1)}$, where $h$ is defined in \eqref{hcall}. Now again using computation like \eqref{sobo-rad-lem-1} we deduce
\begin{align*}
\|f_k\|^2_{D^{1,2}(\m)}=\uv\int_0^{\infty} 
\left( \frac{\partial u }{\partial \varrho}(k\varrho)\right)^2H(\varrho) \varrho^{N-1}\drho=k^{2-N}\uv\int_0^{\infty} \left( \frac{\partial u }{\partial s}(s)\right)^2 H\left(\frac{s}{k}\right)s^{N-1}\ds,
\end{align*}
and then using \eqref{h-lomit} and applying the dominated convergence theorem, we have 
\begin{align}\label{f-3}
    \lim_{k\to \infty}k^{N-2} \|f_k\|^2_{D^{1,2}(\m)} =\|u\|_{D^{1,2}(\rn)}^2.
\end{align}

Now for any $\lambda>0$, we also notice that 
\begin{align}\label{f-4}
    \mu_{u_k}(\lambda)=\int_{\{y\in \rn\, :\, |u(ky)|>\lambda\}}\,{\rm d}y=k^{-N}\int_{\{z\in \rn\, :\, |u(z)|>\lambda\}}\,{\rm d}z=k^{-N}\mu_u(\lambda).
\end{align}

Then by definition $u_k^\star(\mu_{u_k}(\lambda))=\lambda$ and so by substituting $s=k^{-N}\mu_u(\lambda)$ and exploiting \eqref{f-4}, we have 
\begin{align*}
    u_k^\star(s)= \lambda= u^\star(\mu_u(\lambda))=u^\star(k^{N}\, s).
\end{align*}
Using this and making change of variable $t=k^{N} s$, we obtain
\begin{align*}
  \|u_k\|^2_{L^{2^\star,2}(\rn)}=  \int_0^{\infty} \left(s^{\frac{1}{2^\star}}u_k^\star(s)\right)^2\frac{\ds}{s}=k^{2-N}\int_0^{\infty} \left(t^{\frac{1}{2^\star}}u^\star(t)\right)^2\frac{\dt}{t}=k^{2-N}\|u\|^2_{L^{2^\star,2}(\rn)}.
\end{align*}
Now using \eqref{iso-lem}, we have $\|f_k\|^2_{L^{2^\star,2}(\m)}= \|u_k\|^2_{L^{2^\star,2}(\rn)}$ and by the help of above computation and \eqref{f-3}, we deduce
\begin{align*}
    \lim_{k \to \infty}\frac{\|f_k\|^2_{L^{2^\star,2}(\m)}}{\|f_k\|^2_{D^{1,2}(\m)}}= \lim_{k \to \infty}\frac{k^{N-2}\|u_k\|^2_{L^{2^\star,2}(\rn)}}{k^{N-2}\|f_k\|^2_{D^{1,2}(\m)}}= \lim_{k \to \infty}\frac{\|u\|^2_{L^{2^\star,2}(\rn)}}{k^{N-2}\|f_k\|^2_{D^{1,2}(\m)}}=\frac{\|u\|^2_{_{L^{2^\star,2}(\rn)}}}{\|u\|^2_{D^{1,2}(\rn)}},
\end{align*}
which is essentially \eqref{fnl-claim-petri} and this concludes \eqref{f-1}.

In the final part, we will show that the constant $S_{N,2^\star,2}$ is never attained. In fact, if it is attained for some radial non-decreasing function, say $f\in D^{1,2}(\m)$, then taking $\scal(f)\in D^{1,2}(\rn)$ and thanks to \eqref{eq-sobo-rad-lem} and \eqref{iso-lem}, it will give that $S_{N,2^\star,2}$ is attained in the Euclidean case, and this is not possible (refer \cite[Remark~2]{crt}). This completes the thesis.
\end{proof}

\medskip 

\section{Interplay between Cartan–Hadamard and Euclidean Hardy Inequalities}\label{sec-hce}

In this section, we show an equivalence between the geometric and Euclidean frameworks by demonstrating that the Hardy inequality on a Cartan–Hadamard model manifold can be both derived from and sharpened using a suitably weighted Euclidean Hardy inequality. More precisely, we show that the Hardy \emph{deficit} on the manifold admits a quantitative lower bound in terms of a weighted Hardy deficit in the Euclidean space, evaluated for a function naturally associated with the manifold through the Jacobian transformation arising from the radial change of variables  \eqref{jacobian}. This result reveals how the underlying curvature and volume growth characteristics of the model manifold are reflected analytically in the weighted Euclidean structure.

\medskip 

To do this, we first prove a weighted Euclidean Hardy inequality.

\begin{proposition}\label{wg-ecld-hyardy}
Let $N\geq 3$. Suppose $w(x)$ is a non-negative, radially non-decreasing function. Then for every $u\in D^{1,2}(\mathbb{R}^N\setminus \{0\})$ the following holds
    \begin{equation}\label{weighted Euclidean}
         \left(\frac{N-2}{2}\right)^2\int_{\mathbb{R}^N}\frac{|u(x)|^2}{|x|^2}w(x)\dx\leq \int_{\mathbb{R}^N}|\nabla u(x)|^2 w(x)\dx.
    \end{equation}
\end{proposition}
\begin{proof}
We begin by considering radial functions \(u\). It is sufficient to establish the result for any
\[
u \in C_{c,\mathrm{rad}}^{\infty}(\mathbb{R}^N \setminus \{0\}),
\]
since this space is dense in 
\[
D^{1,2}_{\mathrm{rad}}(\mathbb{R}^N \setminus \{0\}).
\]
Let \(u(x) = U(r)\), where \(r = |x|\) denotes the Euclidean distance from the origin. We define 
\begin{align*}
    \phi(r):=(U(r))^2r^{N-2}w(r) \qquad \text{ for all } r>0.
\end{align*}

Notice that $w$ is almost everywhere differentiable. Then, exploiting the compact support of $u$ and integration by parts, we deduce
\begin{align*}
    (N-2)\int_0^{\infty} (U(r))^2 r^{N-3}w(r)\dr=-2\int_0^\infty U(r)U'(r)r^{N-2}w(r)\dr-\int_0^{\infty} (U(r))^2 r^{N-2}w'(r)\dr.
\end{align*}

Then, using $w'\geq 0$ and H\"older inequality in the above, we deduce
\begin{align*}
        (N-2)\int_0^{\infty} (U(r))^2 r^{N-3}w(r)\dr
        \leq &-2\int_0^\infty U(r)U'(r)r^{N-2}w(r)\dr\\
        \leq &2\left(\int_0^\infty (U'(r))^2r^{N-1}w(r)\dr\right)^{\frac{1}{2}}\left(\int_0^{\infty} (U(r))^2r^{N-3}w(r)\dr\right)^{\frac{1}{2}}.
\end{align*}
Therefore, adjusting the terms, we obtain
\begin{equation*}
\left(\frac{N-2}{2}\right)^2 \int_0^{\infty} (U(r))^2r^{N-3}w(r)\dr\leq \int_0^\infty (U'(r))^2r^{N-1}w(r)\dr,
\end{equation*}
which is the required form of \eqref{weighted Euclidean}. For a non-radial, $u\in C^{\infty}_c(\mathbb{R}^N\setminus\{0\})$  we define,
\begin{equation*}
      \bar{u}(|x|)=\left(\frac{1}{\omega_N}\int_{\sn}|u(r,\theta)|^2\dsn\right)^{1/2},
\end{equation*}
 then
\begin{align*}
    \int_{\mathbb{R}^N}\frac{|\bar{u}(|x|)|^2}{|x|^2}w(x)\dx
    =&\uv \int_0^{\infty} |\bar{u}(r)|^2 \frac{w(r)}{r^2}r^{N-1}\dr\\
    =& \uv \int_0^{\infty} \left(\frac{1}{\omega_N}\int_{\sn}|u(r,\theta)|^2\dsn\right)\frac{w(r)}{r^2}r^{N-1}\dr
    = \int_{\mathbb{R}^N} \frac{|u(x)|^2}{|x|^2}w(x)\dx.
\end{align*}
 Also we have
 \begin{equation}\label{radialization comparison}
     \int_{\mathbb{R}^N} |\nabla\bar{u}(x)|^2w(x)\dx\leq \int_{\mathbb{R}^N} |\nabla u(x)|^2w(x)\dx.
 \end{equation}
 For the proof of \eqref{radialization comparison} see, \cite[Lemma~4.1]{pems}. Hence, from the following chain of inequalities
\begin{align*}
  \left(\frac{N-2}{2}\right)^2  \int_{\mathbb{R}^N} \frac{|u(x)|^2}{|x|^2}w(x)\dx= \left(\frac{N-2}{2}\right)^2\int_{\mathbb{R}^N}\frac{|\bar{u}(|x|)|^2}{|x|^2}w(x)\dx & \leq \int_{\mathbb{R}^N} |\nabla\bar{u}(x)|^2w(x)\dx \\
  & \leq \int_{\mathbb{R}^N} |\nabla u(x)|^2w(x)\dx,
\end{align*}
we get our result \eqref{weighted Euclidean} for the compactly supported non-radial functions and thereby for all functions of $D^{1,2}(\mathbb{R}^N\setminus \{0\})$ by a density argument.   
\end{proof}

We are now in a position to present the first main result of this section. The theorem below provides a bridge between the geometric and Euclidean settings by showing that the Hardy inequality on a Cartan–Hadamard model manifold can be derived and strengthened through an appropriately weighted Euclidean Hardy inequality.

\medskip 

\begin{theorem}\label{manifold-eu}
Let $(M,g)$ be a Cartan--Hadamard model manifold with dimension $N\geq 3$. Then for any $u\in D^{1,2}(\m)$ there exists a decreasing function $F\in D^{1,2}_{\text{rad}}(\rn\setminus\{0\})$ such that
\begin{multline*}
     \int_{\m}|\gradg u|^2\dvg-\left(\frac{N-2}{2}\right)^2\int_{\m}\frac{u^2}{r^2}\dvg\\ \geq \int_{\mathbb{R}^N}|\nabla F(x)|^2 (w(x))^{N-1} \dx -\left(\frac{N-2}{2}\right)^2\int_{\mathbb{R}^N}\frac{|F(x)|^2}{|x|^2}(w(x))^{N-1}\dx,
\end{multline*}
   where $w(x)=\frac{\psi(|x|)}{|x|}$ for all $ x\in \rn\setminus\{0\}$.
\end{theorem}

\begin{proof}
Using the Hardy-Littlewood  and  Pólya–Szegő inequalities, we have
\begin{equation}\label{improvement}
    \int_{\m}|\gradg u|^2\dvg-\left(\frac{N-2}{2}\right)^2\int_{\m}\frac{u^2}{r^2}\dvg\geq\int_{\m}|\gradg u^\sharp|^2\dvg-\left(\frac{N-2}{2}\right)^2\int_{\m}\frac{{(u^\sharp)}^2}{r^2}\dvg
\end{equation}
Then, similar to \eqref{1D Etype}, we deduce
\begin{multline}\label{equivalence}
    \int_{\m}|\gradg u^\sharp|^2\dv-\left(\frac{N-2}{2}\right)^2\int_{\m}\frac{{(u^\sharp)}^2}{r^2}\dvg
  \\ = \int_0^{\infty}\left(-\phi'(\varrho)\right)^2\varrho^{\frac{2}{N'}}J(\varrho)\drho-\left(\frac{1}{2^\star}\right)^2\int_0^{\infty} |\phi(\varrho)|^2\varrho^{-\frac{2}{N}}J(\varrho)\drho,
\end{multline}
where $J(\varrho)$ is defined in \eqref{jacobian}, $G(s)=\left(\frac{N}{\uv}\varrho\right)^{\frac{1}{N}}$, $s=V(r)$ and $\phi(\varrho)=(u^\sharp)^\star (s)$.

Recall that $\phi$ is a one-dimensional, non-negative, decreasing, right continuous function. Now we define the $N$-dimensional Euclidean function related to our starting $u$ by
\begin{equation*}
    F(x):= \phi\left(\frac{\uv}{N}|x|^N\right)  \qquad \text{ for all } \quad x\in \rn\setminus\{0\}.
\end{equation*}

Now, we claim that
\begin{equation}\label{claim}
    F^\star_E(s)=\phi(s)\qquad \text{ for all } s>0,
\end{equation}
where $F^\star_E$ is the Euclidean one-dimensional decreasing rearrangement of $F$.

Due to the fact that $\phi_E^\star=\phi$, to prove our claim \eqref{claim}, it is enough to show that $F$ and $\phi$ are equidistributed. Now we define $T: \mathbb{R}^N\to [0,\infty)$ by 
\begin{align*}
  T(x):=\frac{\uv}{N}|x|^N  \qquad \text{ for all } \quad x\in \rn.
\end{align*}
Then $T$ is a measure preserving transformation (see \cite[Example.~1.23]{bnstn}) which implies
\begin{equation}\label{measure equality}
    \mathcal{L}^N\left(T^{-1}(A)\right)= \mathcal{L}^1(A),
\end{equation}
for any measurable $A\subset [0,\infty)$ and $\mathcal{L}^N(\cdot)$ is $N$-dimensional Euclidean Lebesgue measures respectively.

Again we denote
\begin{equation*}
E_F(\lambda)=\left\{x\in \mathbb{R}^N: F(x)>\lambda\right\},\qquad \text{ and }\quad  E_{\phi}(\lambda)=\left\{s> 0: \phi(s)>\lambda\right\}.
\end{equation*}
Now, since $F=\phi \circ T$, therefore, we have $E_F(\lambda)= T^{-1}E_{\phi}(\lambda)$. Hence, thanks to \eqref{measure equality}, we get that $F$ and $\phi$ are equidistributed, and \eqref{claim} follows.

Next, observing the radially symmetric and decreasing property of $F$ along with
\begin{align*}
    F(x)=F^\sharp(x)= F^\star_E\left(\frac{\uv}{N}|x|^N\right)=\phi(\varrho),
\end{align*}
and performing a change of variable $\varrho=\frac{\uv}{N}|x|^N$ in \eqref{equivalence}, we deduce
\begin{align}\label{Euclidean 1D to N dimension}
    &\int_0^{\infty}\left(-\phi'(\varrho)\right)^2\varrho^{\frac{2}{N'}}J(\varrho)\drho-\left(\frac{1}{2^\star}\right)^2\int_0^{\infty} |\phi(\varrho)|^2\varrho^{-\frac{2}{N}}J(\varrho)\drho\nonumber\\&
    =\int_{\mathbb{R}^N}|\nabla F(x)|^2w(x)\dx -\left(\frac{N-2}{2}\right)^2\int_{\mathbb{R}^N}\frac{|F(x)|^2}{|x|^2}w(x)\dx.
\end{align}
Finally, using $\psi(t)\geq t$ for $t>0$ from \eqref{sturm}, we observe $F\in D^{1,2}_{\text{rad}}(\rn\setminus\{0\})$ and using \eqref{Euclidean 1D to N dimension} in \eqref{equivalence} and \eqref{improvement}, and we deduce our result.
\end{proof}

\medskip 

We now turn to prove the Euclidean Hardy inequality obtained from its counterpart on a Cartan–Hadamard manifold. Recall that in the reverse direction, the weight appearing in Theorem~\ref{manifold-eu} was the increasing function $\frac{\psi(t)}{t}$. Hence, when transferring the inequality from the manifold setting back to the Euclidean one, it is natural to anticipate the emergence of a weight with the opposite monotonicity. As we shall see, this is indeed the case. Motivated by this observation, we introduce the reciprocal of the previous weight function, which will play a central role in our subsequent arguments:
\begin{align*}
   \Psi(t):=\frac{t}{\psi(t)} \qquad \text{for all } t>0.
\end{align*}
With this definition in place, we are now ready to establish our next result, which presents a weighted Hardy inequality on Cartan--Hadamard manifolds.

\begin{proposition}\label{prop-mani-euc}
Let $(\m,g)$ be an  $N$-dimensional Cartan–Hadamard model manifold with $N\geq 3$. Then 
\begin{equation}\label{prop-mani-euc-eq}
\left(\frac{N-2}{2}\right)^2\int_{\m} \frac{|u|^2}{r^2}\left(\Psi(r)\right)^{N-1}\dvg\leq \int_{\m}|\gradg u|^2\left(\Psi(r)\right)^{N-1}\dvg
\end{equation}
 holds for every radially decreasing function  $u$ whenever $\lim_{r\to+\infty}u(r)=0$ and r.h.s. of \eqref{prop-mani-euc-eq} is finite.
\end{proposition}
\begin{proof}
Following Proposition~\ref{1D hardy-prop}, the inequality \eqref{prop-mani-euc-eq} can be reformulated in its one-dimensional form. Hence, it is sufficient to establish the inequality
 \begin{equation}\label{mani-euc-eq}
       \left(\frac{N-2}{2}\right)^2\int_0^{\infty}u^\star(s)^2 (G(s))^{N-3}(\psi(G(s))^{1-N}\ds\leq \uv^2\int_0^{\infty} |{u^\star}'(s)|^2\left(\psi \left(G(s)\right)\right)^{(N-1)}(G(s))^{N-1}\ds,
   \end{equation}
   where $r=G(s)$ defined in \eqref{inv-vol} and $u^\star$ is a dimensional rearrangement of $u$.

   For notational simplicity, we will use the following
    \begin{align*}
         p(s):= \left(G(s)\psi (G(s))\right)^{N-1}, \qquad q(s):= (G(s))^{N-3} (\psi (G(s)))^{1-N},
    \end{align*}
    and we define 
    \begin{align*}
        Q(s) = \int_0^s q(\tau)\:{\rm d}\tau \qquad \text{ for all } s>0.
    \end{align*}
    
    To compute the above integral, we make a change of variable $t=G(\tau)$. Then, using the definition of $G$ and making the chain rule, we deduce that for any $s>0$, there holds
    \begin{align*}
        Q(s) = \int_0^s q(\tau)\:{\rm d}\tau = \uv \int_0^{G(s)}t^{N-3}\dt = \frac{\uv}{N-2}(G(s))^{N-2}.
    \end{align*}
    This entails
    \begin{align*}
           \frac{(Q(s))^2}{p(s)q(s)}=\frac{\uv ^2}{(N-2)^2} \qquad \text{ for all } \quad s>0.
    \end{align*}

    Next, the hypothesis $\lim_{r\to+\infty}u(r)=0$ for radially decreasing function $u$ ensures that $u^\star$ which is non-increasing, locally absolutely continuous will also satisfy $\lim_{s\to+\infty}u^\star(s)=0$. Therefore, for every $s>0$, we write
    \begin{align*}
        u^\star(s)^2=-2\int_s^{\infty} u^\star(m){u^\star}'(m)\dm.
    \end{align*}
    Now, multiplying above by $q(s)$ and doing integration w.r.t. $s$ over $(0,\infty)$ and changing the order of integration by Fubini's theorem, we obtain
    \begin{align*}
        \int_0^{\infty} (u^\star(s))^2 q(s)\ds&= -2\int_0^{\infty} \left(\int_s^{\infty}u^\star(m){u^\star}'(m)\dm\right)q(s)\ds \\
        &= 2\int_0^{\infty}u^\star(m)(-{u^\star}'(m))\left(\int_0^m q(s)\ds\right)\dm\\
        &= 2\int_0^{\infty} u^\star(m)(-{u^\star}'(m))Q(m)\dm.
    \end{align*}

    Now applying Young's inequality $2ab\leq \epsilon a^2+\frac{b^2}{\epsilon}$ for $\epsilon>0$, with $a=\sqrt{p(m)}(-{u^\star}'(m))$ and $b= \frac{u^\star(m)Q(m)}{\sqrt{p(m)}}$, we get    
\begin{align*}
    \int_0^{\infty} (u^\star(s))^2 q(s)\ds&\leq\epsilon \int_0^{\infty}(-{u^\star}'(m))^2p(m)\dm+\frac{1}{\epsilon}\int_0^{\infty}\frac{(Q(m) u^\star(m))^2}{p(m)}\dm\\
    &=\epsilon \int_0^{\infty}(-{u^\star}'(m))^2p(m)\dm+\frac{1}{\epsilon}\int_0^{\infty} q(m)(u^\star(m))^2\frac{(Q(m))^2}{p(m)q(m)}\dm\\
    &=\epsilon \int_0^{\infty}(-{u^\star}'(m))^2p(m)\dm+\frac{\uv ^2}{(N-2)^2}\frac{1}{\epsilon}\int_0^{\infty}q(m)(u^\star(m))^2\dm.
\end{align*}
Finally, changing the variable $m$ to $s$ and choosing $\epsilon= \frac{2\uv ^2}{(N-2)^2}$ and arranging terms, we obtain \eqref{mani-euc-eq} and hence \eqref{prop-mani-euc-eq}.
\end{proof}

\medskip 

\begin{theorem}\label{eu-manifold}
    Let $(\m,g)$ be a Cartan--Hadamard model manifold with dimension $N\geq 3$. Then for any $u\in D^{1,2}(\rn\setminus \{0\})$ there exists a radially decreasing function $F$ defined on $\m\setminus\{x_0\}$ with limit zero at infinity, and the following holds
    \begin{multline*}
         \int_{\mathbb{R}^N}|\nabla u|^2 \dx-\left(\frac{N-2}{2}\right)^2\int_{\mathbb{R}^N}\frac{|u|^2}{|x|^2}\dx \\ \geq  \int_{\m}|\gradg F|^2\left(\Psi(r)\right)^{N-1}\dvg-\left(\frac{N-2}{2}\right)^2\int_{\m} \frac{|F|^2}{r^2}\left(\Psi(r)\right)^{N-1}\dvg.
    \end{multline*}
\end{theorem}

\begin{proof}
    Start with $u\in D^{1,2}(\rn\setminus \{0\})$. Using the Hardy-Littlewood  and  Pólya–Szegő inequality in the l.h.s. of \eqref{eu}, and taking the transformation $\varrho= \frac{\uv}{N}|x|^N$, $\phi(\varrho)= (u^\sharp)^\star(\varrho)$, it reduces to the following well-known one-dimensional form
\begin{align}\label{eu-1}
    \int_{\mathbb{R}^N}|\nabla u|^2 \dx-\left(\frac{N-2}{2}\right)^2\int_{\mathbb{R}^N}\frac{|u|^2}{|x|^2}\dx\geq \int_0^{\infty}\left(-\phi'(\varrho)\right)^2\varrho^{\frac{2}{N'}}\drho-\left(\frac{1}{2^\star}\right)^2\int_0^{\infty} |\phi(\varrho)|^2\varrho^{-\frac{2}{N}}\drho,
\end{align}
where $\sharp,\,\star$ denote the $N$-dimensional Euclidean symmetric decreasing rearrangement and the one-dimensional Euclidean decreasing rearrangement, respectively. Now, recall the function $G$ defined in \eqref{inv-vol} , and for $s>0$, we define the function
\begin{align*}
    w(s):=\phi(\varrho) \qquad \text{ where } \quad G(s)=\left(\frac{N}{\uv}\varrho\right)^{\frac{1}{N}}.
\end{align*}

Now we are ready to define our required function on the Cartan--Hadamard model manifold $\m$, associated to our starting $u$ below
\begin{equation*}
    F(x):= w\left(V(r)\right)  \qquad \text{ for all } \quad x\in \m \setminus\{x_0\} \quad \text{ and } \quad r=\text{dist}(x,x_0).
\end{equation*}
Now, it is straightforward to see that the r.h.s. of \eqref{eu-1} reduces to
\begin{align}\label{eu-2}
   \ \uv^2\int_0^{\infty} |w'(s)|^2\left(\psi \left(G(s)\right)\right)^{2(N-1)}\left(\Psi(G(s))\right)^{N-1}\ds-  \left(\frac{N-2}{2}\right)^2\int_0^{\infty}\frac{w(s)^2}{\left(G(s)\right)^2}\left(\Psi(G(s))\right)^{N-1}\ds. 
\end{align}

Now, similarly to the proof of Theorem~\ref{manifold-eu}, we can show $F^\star(s)=w(s)$ for every $s>0$, and here $\star$ denotes one-dimensional decreasing rearrangement on the Cartan--Hadamard model manifold. Therefore, the r.h.s. of \eqref{eu} becomes \eqref{eu-2} and the proof follows. 
\end{proof}

\begin{remark}
{\rm 
From Proposition~\ref{prop-mani-euc}, we observe that the right-hand side of \eqref{eu} is nonnegative. This yields an improvement of the classical Euclidean Hardy inequality, now strengthened by a weighted Hardy inequality adapted to the geometry of a general Cartan–Hadamard model manifold $(\m,g)$.
}
\end{remark}

\medskip 


\section*{Acknowledgments} 
A.~Banerjee is supported by the Doctoral Fellowship of the Indian Statistical Institute, Delhi Centre. The research of D.~Ganguly is partially supported
by the ANRF MATRICS Research Grant (MTR/2023/000331) and ANRF-Advanced Research
Grant (ANRF/ARG/2025/000572/MS). P.~Roychowdhury is supported by the Indian Institute of Technology Hyderabad through the Seed Grant under research grant No. SG/IITH/F391/2026-27/SG-230.

\section*{Competing interests} 
The authors have no competing interests to declare that are relevant
to the content of this article.

\section*{Data availability statement}
Data sharing is not applicable to this article as no data sets were generated or analysed during the current study.

\section*{AI statement}
The strategy and arguments are entirely human-generated, and the paper does not contain AI-generated text. The authors take full responsibility for the entire content of the paper. AI-based tools were used only for language-related assistance, such as spell-checking and minor grammatical corrections, and not for developing the mathematical arguments, results, or conclusions.

\end{document}